\documentclass[11pt]{article}
\usepackage{latexsym}
\usepackage{amsfonts,amssymb,amsmath}

\usepackage[yyyymmdd]{datetime}

\newtheorem{thm}{Theorem}[section]
\newtheorem{cor}[thm]{Corollary}
\newtheorem{lem}[thm]{Lemma}
\newtheorem{pro}[thm]{Proposition}
\newtheorem{defn}[thm]{Definition}

\title{Monoid generalizations of the Richard Thompson groups}

\author{ J.C.\ Birget
   \thanks{\, Earlier versions of this paper appeared in 
    http://arxiv.org/abs/0704.0189 (v1 in April 2007, and v2 in April 2008),
    and also appeared in reference \cite{BiMonTh}. Supported by NSF grant
    CCR-0310793.}
       }

\date{\footnotesize{24 Jan.\ 2016}}

\begin{document}
\maketitle

\begin{abstract}
The groups $G_{k,1}$ of Richard Thompson and Graham Higman can be generalized
in a natural way to monoids, that we call $M_{k,1}$, and to inverse monoids,
called ${\it Inv}_{k,1}$; this is done by simply generalizing bijections to
partial functions or partial injective functions.
The monoids $M_{k,1}$ have connections with circuit complexity (studied in 
another paper). 
Here we prove that $M_{k,1}$ and ${\it Inv}_{k,1}$ are congruence-simple for 
all $k$. Their Green relations $\cal J$ and $\cal D$ are characterized: 
$M_{k,1}$ and ${\it Inv}_{k,1}$ are $\cal J$-0-simple, and they have $k-1$ 
non-zero $\cal D$-classes. 
They are submonoids of the multiplicative part of the Cuntz algebra 
${\cal O}_k$. 
They are finitely generated, and their word problem over any finite 
generating set is in {\sf P}. Their word problem is {\sf coNP}-complete over 
certain infinite generating sets.\footnote{\, {\bf Changes in this version:} 
Section 4 has been thoroughly revised, and errors have been corrected;
however, the main results of Section 4 do not change.
The main changes are in Theorem 4.5, Definition 4.5A (the concept of a 
{\em normal} right-ideal morphism), and the final proof of Theorem 4.13.  
Sections 1, 2, and 3 are unchanged, except for the proof of Theorem 2.3, which was incomplete; a complete proof was published in the Appendix of 
reference \cite{BiBern}, and is also given here.}
\end{abstract}


\section{Thompson-Higman monoids}

Since their introduction by Richard J.\ Thompson in the mid 1960s 
\cite{Th0, McKTh, Th}, the Thompson groups have had a great impact on 
infinite group theory. Graham Higman generalized the Thompson groups to 
an infinite family \cite{Hig74}. These groups and some of their subgroups 
have appeared in many contexts and have been widely studied; see for example 
\cite{CFP, BrinSqu, Dehornoy, BrownGeo, GhysSergiescu, GubaSapir, Brin97, 
BCST, LawsonPolycyclic}.

The definition of the Thompson-Higman groups lends itself easily to 
generalizations to inverse monoids and to more general monoids.
These monoids are also generalizations of the finite symmetric monoids 
(of all functions on a set), and this leads to connections with 
circuit complexity; more details on this appear in 
\cite{BiDistor, BiFact, BiCoNP}.

By definition the Thompson-Higman group $G_{k,1}$ consists of all maximally
extended isomorphisms between finitely generated essential right
ideals of $A^*$, where $A$ is an alphabet of cardinality $k$.
The multiplication is defined to be composition followed by maximal extension:
for any $\varphi, \psi \in G_{k,1}$, we have  \ 
$\varphi \cdot \psi = $ max$(\varphi \circ \psi)$.
Every element $\varphi \in G_{k,1}$ can also be given by a bijection 
$\varphi: P \to Q$ where $P, Q \subset A^*$ are two finite maximal prefix 
codes over $A$; this bijection can be described concretely by a finite 
function {\it table}. For a detailed definition according to this approach,
see \cite{BiThomps} (which is also similar to \cite{Scott}, but with a
different terminology); moreover, Subsection 1.1 gives all the needed 
definitions. 

It is natural to generalize the maximally extended {\em isomorphisms} between
finitely generated essential right ideals of $A^*$  to {\em homomorphisms}, 
and to drop the requirement that the right ideals be essential. It will turn 
out that this generalization leads to interesting monoids, or inverse monoids, 
which we call Thompson-Higman monoids. Our generalization of the
Thompson-Higman groups to monoids will also generalize the embedding of these
groups into the Cuntz algebras \cite{BiThomps, Nekrash}, which provides
an additional motivation for our definition.
Moreover, since these homomorphisms are close to being arbitrary finite
string transformations, there is a connection between these monoids 
and combinational boolean circuits; the study of the connection 
between Thompson-Higman groups and circuits was started in 
\cite{BiCoNP, BiFact} and will be developed more generally for monoids 
in \cite{BiDistor}; the present paper lays some of the foundations for
\cite{BiDistor}.


\subsection{Definition of the Thompson-Higman groups and monoids}

Before defining the Thompson-Higman monoids we need some basic definitions,
that are similar to the introductory material that is needed for defining 
the Thompson-Higman groups $G_{k,1}$; we follow \cite{BiThomps} (which is 
similar to \cite{Scott}).  
We use an alphabet $A$ of cardinality $|A| = k$, and we list its elements 
as $A = \{a_1, \ldots, a_k\}$. Let $A^*$ denote the set of all finite
{\it words} over $A$ (i.e., all finite sequences of elements of $A$); this 
includes the {\it empty word} $\varepsilon$. The {\it length} of $w \in A^*$ 
is denoted by $|w|$; let $A^n$ denote the set of words of length $n$. For 
two words $u,v \in A^*$ we denote their {\it concatenation} by $uv$ or by
$u \cdot v$; for sets $B, C \subseteq A^*$ the concatenation is
$BC = \{uv : u \in B, v \in C\}$.
A {\it right ideal} of $A^*$ is a subset $R \subseteq A^*$ such that
$RA^* \subseteq R$. A generating set of a right ideal $R$ is a set $C$ such 
that $R$ is the intersection of all right ideals that contain $C$;
equivalently, $R = CA^*$. A right ideal $R$ is called {\it essential} iff 
$R$ has a non-empty intersections with every right ideal of $A^*$.
For words $u,v \in A^*$, we say that $u$ is a {\it prefix} of $v$ iff there
exists $z \in A^*$ such that $uz = v$. A {\it prefix code} is a subset
$C \subseteq A^*$ such that no element of $C$ is a prefix of another element
of $C$. A prefix code is {\it maximal} iff it is not a strict subset of
another prefix code.
One can prove that a right ideal $R$ has a unique minimal (under inclusion)
generating set, and that this minimal generating set is a prefix code; this
prefix code is maximal iff $R$ is an essential right ideal.

For right ideals $R' \subseteq R \subseteq A^*$ we say that $R'$ is 
{\em essential in} $R$ iff  $R'$ intersects all right subideals of $R$ in a 
non-empty way.

\medskip

{\em Tree interpretation:}  The free monoid $A^*$ can be pictured by its 
right Cayley graph, which is the rooted infinite regular $k$-ary tree with 
vertex set $A^*$ and edge set 
$\{(v,va): v \in A^*, a \in A\}$. We simply call this the {\em tree of} 
$A^*$. It is a directed tree, with all paths moving away from the root 
$\varepsilon$ (the empty word); by ``path'' we will always mean a directed 
path. A word $v$ is a prefix of a word $w$ iff $v$ is is an ancestor of $w$ 
in the tree. 
A set $P$ is a prefix code iff no two elements of $P$ are on the same path.   
A set $R$ is a right ideal iff any path that starts in $R$ has all its 
vertices in $R$.
The prefix code that generates $R$ consists of the elements of $R$ that are
maximal (within $R$) in the prefix order, i.e., closest to the root 
$\varepsilon$. 
A finitely generated right ideal $R$ is essential iff every infinite path 
of the tree eventually reaches $R$ (and then stays in it from there on). 
Similarly, a finite prefix code $P$ is maximal iff any infinite path starting 
at the root eventually intersects $P$. 
For two finitely generated right ideals $R' \subset R$, $R'$ is essential in
$R$ iff any infinite path starting in $R$ eventually reaches  $R'$ (and then
stays in $R'$ from there on). In other words for finitely generated right 
ideals $R' \subseteq R$, $R'$ is essential in $R$ iff
$R'$ and $R$ have the same ``ends''. 
For the prefix tree of $A^*$ we can consider also the ``boundary'' 
$A^{\omega}$ (i.e., all infinite words), a.k.a.\ the {\em ends} of the tree.
In Thompson's original definition
\cite{Th0,Th}, $G_{2,1}$ was given by a total action on $\{0,1\}^{\omega}$.
In \cite{BiThomps} this total action was extended to a partial action on
$A^* \cup A^{\omega}$; the partial action on $A^* \cup A^{\omega}$ is
uniquely determined by the total action on $A^{\omega}$; it is also uniquely
determined by the partial action on $A^*$. Here, as in \cite{BiThomps}, we
only use the partial action on $A^*$.

\begin{defn}  \label{homo_right_ideals} \ 
A {\em right ideal homomorphism} of $A^*$ is a total function 
$\varphi: R_1 \to A^*$ such that $R_1$ is a right ideal of $A^*$, and for 
all $x_1 \in R_1$ and all $w \in A^*$:
 \ $\varphi(x_1w) = \varphi(x_1) \ w$.
\end{defn}
For any partial function $f:A^* \to A^*$, let Dom$(f)$ denote the domain and 
let Im$(f)$ denote the image (range) of $f$.
For a right ideal homomorphism $\varphi: R_1 \to A^*$ it is easy to see 
that the image ${\rm Im}(\varphi)$ is also right ideal of $A^*$, which is 
finitely generated (as a right ideal) if the domain 
$R_1 = {\rm Dom}(\varphi)$ is finitely generated. 

A right ideal homomorphism $\varphi: R_1 \to R_2$, where
$R_1 = {\rm Dom}(\varphi)$ and $R_2 = {\rm Im}(\varphi)$, can be described 
by a total surjective function $P_1 \to S_2$, with $P_1, S_2 \subset A^*$; 
here $P_1$ is the 
prefix code (not necessarily maximal) that generates $R_1$ as a right ideal, 
and $S_2$ is a set (not necessarily a prefix code) that generates $R_2$
as a right ideal; so $R_1 = P_1 A^*$ and $R_2 = S_2 A^*$.  
The function $P_1 \to S_2$ corresponding to
$\varphi: R_1 \to R_2$ is called the {\em table} of $\varphi$. The prefix 
code $P_1$ is called the {\em domain code} of $\varphi$ and we write 
$P_1 = {\rm domC}(\varphi)$.  When $S_2$ is a prefix code we call $S_2$ the 
{\em image code} of $\varphi$ and we write $S_2 = {\rm imC}(\varphi)$.
We denote the {\em table size} of $\varphi$ (i.e., the cardinality of
${\rm domC}(\varphi)$) by $\|\varphi\|$.

\begin{defn}  \label{special_homo_right_ideals} \ 
An injective right ideal homomorphism is called a right ideal 
{\em isomorphism}.
A right ideal homomorphism $\varphi: R_1 \to R_2$ is called {\em total} iff
the domain right ideal $R_1$ is essential. And $\varphi$ is called 
{\em surjective} iff the image right ideal $R_2$ is essential.
\end{defn}
The table $P_1 \to P_2$ of a right ideal isomorphism $\varphi$ is a bijection
between prefix codes (that are not necessarily maximal). The table 
$P_1 \to S_2$ of a {\em total} right ideal homomorphism is a function from a 
{\em maximal} prefix code to a set, and the table $P_1 \to S_2$ of a 
surjective right ideal homomorphism is a function from a prefix code to a 
set that generates an essential right ideal.  
The word ``total'' is justified by the fact that if a homomorphism $\varphi$ 
is total (and if ${\rm domC}(\varphi)$ is finite) then $\varphi(w)$ is defined 
for every word that is long enough (e.g., when $|w|$ is longer than the 
longest word in the domain code $P_1$); equivalently, $\varphi$ is defined 
from some point onward on every infinite path in the tree of $A^*$ starting 
at the root.

\begin{defn}  \label{restriction_extension}
An {\em essential restriction} of a right ideal homomorphism 
$\varphi: R_1 \to A^*$ is a right ideal homomorphism $\Phi: R'_1 \to A^*$ 
such that $R'_1$ is essential in $R_1$, and such that for all 
$x'_1 \in R'_1$: \ $\varphi(x'_1) = \Phi(x'_1)$.

We say that $\varphi$ is an {\em essential extension} of $\Phi$ iff $\Phi$ is  
an  essential restriction of $\varphi$.
\end{defn}
Note that if $\Phi$ is an essential restriction of $\varphi$ then 
$R'_2 = {\rm Im}(\Phi)$ will automatically be essential in
$R_2 = {\rm Im}(\varphi)$. Indeed, if $I$ is any non-empty right subideal of 
$R_1$ then $I \cap R'_1 \neq \varnothing$, hence \ 
$\varnothing \neq \Phi(I \cap R'_1)$
$\subseteq \Phi(I) \, \cap \, \Phi(R'_1)$ $= \Phi(I) \cap R'_2$; 
moreover, any right subideal $J$ of $R_2$ is of the form 
$J = \Phi(I)$ where $I = \Phi^{-1}(J)$ is a right subideal of $R_1$; 
hence, for any right subideal $J$ of $R_2$, $\varnothing \neq J \cap R'_2$.  

\begin{pro} \label{extension_formula} \     
{\rm (1)}  Let $\varphi, \Phi$ be homomorphisms between {\em finitely 
generated} right ideals of $A^*$, where $A = \{a_1, \ldots, a_k\}$. 
Then $\Phi$ is an essential restriction of $\varphi$ 
iff $\Phi$ can be obtained from $\varphi$ by starting from the table of
$\varphi$ and applying a finite number of {\em restriction steps} of the 
following form: Replace $(x,y)$ in a table by 
$\{(xa_1,ya_1), \ldots, (xa_k,ya_k)\}$.  \\     
{\rm (2)} Every homomorphism between finitely generated right ideals of 
$A^*$ has a {\em unique maximal essential} extension. 
\end{pro}
{\bf Proof.} \ (1) Consider a homomorphism between finitely generated
right ideals $\varphi: R_1 \to R_2$, let $P_1$ be the finite prefix code
that generates the right ideal $R_1$, and let $S_2 = \varphi(P_1)$, so $S_2$
generates the right ideal $R_2$.

If $x \in P_1$ and $y = \varphi(x) \in S_2$ then (since $\varphi$ is a right
ideal homomorphism), $ya_i = \varphi(xa_i)$ for $i = 1, \ldots, k$. Then
$R_1 - \{x\}$ is a right ideal which is essential in $R_1$, and $R_1 - \{x\}$
is generated by $(P_1 - \{x\}) \, \cup \, \{xa_1, \ldots, xa_k\}$. 
Indeed, in the tree of $A^*$ every downward directed path 
starting at vertex $x$ goes through one of the vertices $xa_i$. Thus, 
removing $(x,y)$ from the graph of $\varphi$ is an essential restriction; 
for the table of $\varphi$, the effect is to replace the entry $(x,y)$ by 
the set of entries $\{(xa_1,ya_1), \ldots, (xa_k,ya_k)\}$.
If finitely many restriction steps of the above type are carried out, the
result is again an essential restriction of $\varphi$.

Conversely, let us show that if $\Phi$ is an essential restriction of 
$\varphi$ then $\Phi$ can be obtained by a finite number of replacement steps 
of the form ``replace $(x,y)$ by $\{(xa_1,ya_1), \ldots, (xa_k,ya_k)\}$ in 
the table''. 

Using the tree of $A^*$ we have: If $R$ and $R'$ are right ideals of $A^*$ 
generated by the finite prefix codes $P$, respectively $P'$, and if $R'$
is essential in $R$ then every infinite path from $P$ intersects $P'$.
It follows from this characterization of essentiality and from the finiteness 
of $P_1$ and $P'_1$ that $R_1 - R'_1$ is finite. Hence $\varphi$ and $\Phi$ 
differ only in finitely many places, i.e., one can transform $\varphi$ into 
$\Phi$ in a finite number of restriction steps. 

So, the restriction $\Phi$ of $\varphi$ is obtained by removing a finite 
number of pairs $(x,y)$ from $\varphi$; however, not every such removal 
leads to a right ideal homomorphism or  an essential restriction of $\varphi$.
If $(x_0, y_0)$ is removed from $\varphi$ then $x_0$ is removed from $R_1$
(since $\varphi$ is a function). Also, since $R'_1$ is a right ideal, when 
$x_0$ is removed then all prefixes of $x_0$ (equivalently, all ancestor
vertices of $x_0$ in the tree of $A^*$) have to be removed. 
So we have the following removal rule (still assuming that domain and image
right ideals are finitely generated):

\smallskip

\noindent 
{\it If $\Phi$ is an essential restriction of $\varphi$ then $\varphi$ can 
be transformed into $\Phi$ by removing a finite set of strings from $R_1$, 
with the following restriction: If a string $x_0$ is removed then all prefixes 
of $x_0$ are also removed from $R_1$; moreover, $x_0$ is removed from $R_1$ 
iff $(x_0, \varphi(x_0))$ is removed from $\varphi$.
}

\smallskip

As a converse of this rule, we claim that if the transformation from 
$\varphi$ to $\Phi$ is done according to this rule, then $\Phi$ is an 
essential restriction of $\varphi$. 
Indeed, $\Phi$ will be a right ideal homomorphism: if $\Phi(x_1)$ 
is defined then $\Phi(x_1z)$ will also be defined (if it were not, the prefix
$x_1$ of $x_1z$ would have been removed), and $\Phi(x_1z) = \varphi(x_1z) = $
$\varphi(x_1) \ z = \Phi(x_1) \ z$.
Moreover, ${\rm Dom}(\Phi) = R'_1$  will be essential in $R_1$: every 
directed path starting at $R_1$ eventually meets $R'_1$ because only
finitely many words were removed from $R_1$ to form $R'_1$. 
Hence by the tree characterization of essentiality,  $R'_1$ is essential in
$R_1$.

In summary, if $\Phi$ is an essential restriction of $\varphi$ then $\Phi$ 
is obtained from $\varphi$ by a finite sequence of steps, each of which removes
one pair $(x,\varphi(x))$. In ${\rm Dom}(\varphi)$ the string $x$ is removed.
The domain code becomes $(P_1 - \{x\}) \, \cup \, \{xa_1, \ldots, xa_k\}$,
since $\{xa_1, \ldots, xa_k\}$ is the set of children of $x$ in the tree of
$A^*$. This means that in the table of $\varphi$, the pair $(x,\varphi(x))$
is replaced by 
$\{(xa_1, \varphi(x) \, a_1), \ldots, (xa_k, \varphi(x) \, a_k)\}$.

\smallskip

\noindent (2) Uniqueness of the maximal essential extension: \  
By (1) above, essential extensions are obtained by the set of rewrite 
rules of the form \ $\{(xa_1,ya_1), \ldots, (xa_k,ya_k)\} \to (x,y)$, 
applied to tables. 
This rewriting system is {\it locally confluent} (because different rules 
have non-overlapping left sides) and {\it terminating} (because they decrease 
the length); hence maximal essential extensions exist and are unique. 
 \ \ \ $\Box$ 

\medskip

\noindent 
Proposition \ref{extension_formula} yields another tree interpretation of 
essential restriction: Assume first that a total 
order $a_1 < a_2 < \ldots < a_k$ has been chosen for the alphabet $A$; this 
means that the tree of $A^*$ is now an {\em oriented} rooted tree, 
i.e., the children of each vertex $v$ have a total order (namely,
$va_1 < va_2 < \ldots < va_k$).  The rule ``replace $(x,y)$ in the table by
$\{(xa_1,ya_1), \ldots, (xa_k,ya_k)\}$'' has the following tree 
interpretation: Replace $x$ and $y = \varphi(x)$ by the children of $x$,
respectively of $y$, matched according to the order of the children. 

\medskip

\noindent {\bf Important remark:} \\   
As we saw, every right ideal homomorphism can be described by a table 
$P \to S$ where $P$ is a prefix code and $S$ is a set. But we also have: 
Every right ideal homomorphism $\varphi$ has an essential restriction 
$\varphi'$ whose table $P' \to Q'$ is such that {\em both $P'$ and $Q'$ are 
prefix codes}; moreover, $Q'$ can be chosen to be a subset of 
$A^n$ for some $n \leq {\rm max}\{ |s| : s \in S\}$.  
  Example (with alphabet $A = \{a,b\}$): \\    
$\left( \hspace{-.08in} \begin{array}{r|r} 
a &  b \\ 
a & aa  
\end{array} \hspace{-.08in} \right) $ \ has an essential restriction \    
$ \left( \hspace{-.08in} \begin{array}{r|r|r} 
aa & ab &  b \\   
aa & ab & aa  
\end{array} \hspace{-.08in} \right) $. 
Theorem 4.5B gives a tighter result with polynomial bounds.

\begin{defn} \label{Tomps_part_funct_monoid} \
The Thompson-Higman {\em partial function monoid} $M_{k,1}$ consists of all
maximal essential extensions of homomorphisms between finitely generated
right ideals of $A^*$. The multiplication is composition followed by 
maximal essential extension.
\end{defn}
In order to prove {\it associativity} of the multiplication of $M_{k,1}$
we define the following and we prove a few Lemmas.

\begin{defn} \label{congruence} \   
By RI$_k$ we denote the monoid of all right ideal homomorphisms between 
finitely generated right ideals of $A^*$, with function composition as 
multiplication. We consider the equivalence relation $\equiv$ defined for 
$\varphi_1, \varphi_2 \in {\it RI}_k$ by: \ $\varphi_1 \equiv \varphi_2$ 
 \ iff \ ${\rm max}(\varphi_1) = {\rm max}(\varphi_2)$.
\end{defn}
It is easy to prove that ${\it RI}_k$ is closed under composition. Moreover,
by existence and uniqueness of the maximal essential extension 
(Prop.\ \ref{extension_formula}(2)) each $\equiv$-equivalence class contains 
exactly one element of $M_{k,1}$. We want to prove:

\begin{pro} \label{congruence_is_congruence} \   
The equivalence relation $\equiv$ is a monoid congruence on ${\it RI}_k$, 
and $M_{k,1}$ is isomorphic (as a monoid) to ${\it RI}_k/\!\equiv$. 
 \ Hence, $M_{k,1}$ is associative.
\end{pro}
First some Lemmas. 

\begin{lem} \label{inters_essent_ideals} \   
If $R'_i \subseteq R_i$ ($i = 1,2$) are finitely generated right ideals with
$R'_i$ essential in $R_i$, then $R'_1 \cap R'_2$ is essential in 
$R_1 \cap R_2$.
\end{lem}
{\bf Proof.} We use the tree characterization of essentiality. Any infinite 
path $p$ in $R_1 \cap R_2$ is also in $R_i$ ($i = 1,2$), hence $p$ eventually
enters into $R'_i$. Thus $p$ eventually meets $R'_1$ and $R'_2$, i.e.,
$p$ meets $R'_1 \cap R'_2$. \ \ \  $\Box$

\begin{lem} \label{dom_equal_im} \ 
All $\varphi_1, \varphi_2 \in {\it RI}_k$ have restrictions 
$\Phi_1, \Phi_2 \in {\it RI}_k$ (not necessarily essential restrictions) 
such that: \\   
$\bullet$ \ \ $\Phi_2 \circ \Phi_1 = \varphi_2 \circ \varphi_1$, \ and \\   
$\bullet$ \ \ ${\rm Dom}(\Phi_2) = {\rm Im}(\Phi_1) = $
     ${\rm Dom}(\varphi_2) \, \cap \, {\rm Im}(\varphi_1)$.  
\end{lem}
{\bf Proof.} Let $R = {\rm Dom}(\varphi_2) \cap {\rm Im}(\varphi_1)$.  This 
is a right ideal which is finitely generated since ${\rm Dom}(\varphi_2)$
and ${\rm Im}(\varphi_1)$ are finitely generated (see Lemma 3.3 of 
\cite{BiThomps}). Now we restrict $\varphi_1$ to $\Phi_1$ in such a way that 
${\rm Im}(\Phi_1) = R$ and ${\rm Dom}(\Phi_1) = \varphi_1^{-1}(R)$, and we 
restrict $\varphi_2$ to $\Phi_2$ in such a way that ${\rm Dom}(\Phi_2) = R$ 
and ${\rm Im}(\Phi_2) = \varphi_2(R)$. Then $\Phi_2 \circ \Phi_1(.)$ 
and $\varphi_2 \circ \varphi_1(.)$ agree on $\varphi_1^{-1}(R)$; moreover,
${\rm Dom}(\Phi_2 \circ \Phi_1) = \varphi_1^{-1}(R)$. Since 
$\varphi_2 \circ \varphi_1(x)$ is only defined when $\varphi_1(x) \in R$,
we have $\Phi_2 \circ \Phi_1 = \varphi_2 \circ \varphi_1$. 
Also, by the definition of $R$ we have 
${\rm Dom}(\Phi_2) = {\rm Im}(\Phi_1)$. 
 \ \ \  $\Box$

\begin{lem} \label{max_max} \   
For all $\varphi_1, \varphi_2 \in {\it RI}_k$ we have: \ \   

\medskip

 \ \ \ \ \ ${\rm max}(\varphi_2 \circ \varphi_1) = $
   ${\rm max}({\rm max}(\varphi_2) \circ \varphi_1) = $
   ${\rm max}(\varphi_2 \circ {\rm max}(\varphi_1))$.
\end{lem}
{\bf Proof.}  We only prove the first equality; the proof of the second one
is similar. By Lemma \ref{dom_equal_im} we can restrict $\varphi_1$ and
$\varphi_2$ to $\varphi'_1$, respectively $\varphi'_2$, so that \  
$\varphi'_2 \circ \varphi'_1 = \varphi_2 \circ \varphi_1$, 
and \ ${\rm Dom}(\varphi'_2) = {\rm Im}(\varphi'_1) = $
 ${\rm Dom}(\varphi_2) \cap {\rm Im}(\varphi_1)$; let 
$R' = {\rm Dom}(\varphi_2) \cap {\rm Im}(\varphi_1)$.

Similarly we can restrict $\varphi_1$  and ${\rm max}(\varphi_2)$ to 
$\varphi''_1$, respectively $\varphi''_2$, so that \  
$\varphi''_2 \circ \varphi''_1 = {\rm max}(\varphi_2) \circ \varphi_1$,
and \ ${\rm Dom}(\varphi''_2) = {\rm Im}(\varphi''_1) = $
 ${\rm Dom}({\rm max}(\varphi_2)) \cap {\rm Im}(\varphi_1)$; let 
$R'' = {\rm Dom}({\rm max}(\varphi_2)) \cap {\rm Im}(\varphi_1)$.

Obviously, $R' \subseteq R''$ (since $\varphi_2$ is a restriction of
${\rm max}(\varphi_2)$). Moreover, $R'$ is essential in $R''$, by Lemma
\ref{inters_essent_ideals}; indeed, ${\rm Dom}(\varphi_2)$ is essential in  
${\rm Dom}({\rm max}(\varphi_2))$ since ${\rm max}(\varphi_2)$ is an essential
extension of $\varphi_2$. 
Since $R'$ is essential in $R''$, $\varphi_2 \circ \varphi_1$ is an essential
restriction of ${\rm max}(\varphi_2) \circ \varphi_1$. Hence by uniqueness 
of the maximal essential extension,  
${\rm max}({\rm max}(\varphi_2) \circ \varphi_1) = $
   ${\rm max}(\varphi_2 \circ {\rm max}(\varphi_1))$.
 \ \ \ $\Box$
  
\medskip 

\noindent {\bf Proof of Prop.\ \ref{congruence_is_congruence}:} \   
If $\varphi_2 \equiv \psi_2$ then, by definition, 
${\rm max}(\varphi_2) = {\rm max}(\psi_2)$, hence by Lemma \ref{max_max}: 

\smallskip

 ${\rm max}(\varphi_2 \circ \varphi) = $
 ${\rm max}({\rm max}(\varphi_2) \circ \varphi) = $
 ${\rm max}({\rm max}(\psi_2) \circ \varphi) = $
${\rm max}(\psi_2 \circ \varphi)$, 

\smallskip

\noindent for all $\varphi \in {\it RI}_k$. Thus (by the definition of 
$\equiv$), $\varphi_2 \circ \varphi \equiv \psi_2 \circ \varphi$, so $\equiv$
is a right congruence.
Similarly one proves that $\equiv$ is a left congruence. Thus, 
${\it RI}_k/\!\!\equiv \,$ is a monoid. 

Since every $\equiv$-equivalence class contains exactly one element of 
$M_{k,1}$ there is a one-to-one correspondence between 
${\it RI}_k/\!\!\equiv \, $
and $M_{k,1}$. Moreover, the map 
 \ $\varphi \in {\it RI}_k \longmapsto {\rm max}(\varphi) \in M_{k,1}$ \ is
a homomorphism, by Lemma \ref{max_max} and by the definition of 
multiplication in $M_{k,1}$.
Hence ${\it RI}_k/\!\!\equiv \, $ is isomorphic to $M_{k,1}$. 
 \ \ \ $\Box$
 

\subsection{Other Thompson-Higman monoids}

We now introduce a few more families of Thompson-Higman monoids, whose
definition comes about naturally in analogy with $M_{k,1}$.

\begin{defn} \label{Tomps_monoid_variants} \
The Thompson-Higman {\em total} function monoid ${\it tot}M_{k,1}$
and the Thompson-Higman {\em surjective} function monoid ${\it sur}M_{k,1}$
consist of maximal essential extensions of homomorphisms between finitely 
generated right ideals of $A^*$ where the domain, respectively, the image 
ideal, is an {\em essential} right ideal.

The Thompson-Higman {\em inverse} monoid ${\it Inv}_{k,1}$ consists of all
maximal essential extensions of isomorphisms between finitely generated
(not necessarily essential) right ideals of $A^*$.
\end{defn}
Every element $\varphi \in {\it tot}M_{k,1}$ can be described by a function
$P \to Q$, called the {\it table} of $\varphi$, where $P, Q \subset A^*$ with 
$P$ a finite {\em maximal} prefix code over $A$. 
A similar description applies to ${\it sur}M_{k,1}$ but now 
with $Q$ a finite maximal prefix code.
Every $\varphi \in {\it Inv}_{k,1}$ can be described by a bijection
$P \to Q$ where $P, Q \subset A^*$ are two finite prefix codes
(not necessarily maximal).

It is easy to prove that essential extension and restriction of right 
ideal homomorphisms, as well as composition of such homomorphisms, preserve 
injectiveness, totality, and surjectiveness. 
Thus ${\it tot}M_{k,1}$, ${\it sur}M_{k,1}$, and ${\it Inv}_{k,1}$ are
submonoids of $M_{k,1}$.

We also consider the intersection ${\it tot}M_{k,1} \cap {\it sur}M_{k,1}$, 
i.e., the monoid of all maximal essential extensions of homomorphisms 
between finitely generated essential right ideals of $A^*$; we denote this
monoid by ${\it totsur}M_{k,1}$.
The monoids $M_{k,1}$, ${\it tot}M_{k,1}$, ${\it sur}M_{k,1}$, and 
${\it totsur}M_{k,1}$ are {\em regular} monoids. (A monoid $M$ is regular iff
for every $m \in M$ there exists $x \in M$ such that $m x m = m$.) 
The monoid ${\it Inv}_{k,1}$ is an inverse monoid. (A monoid $M$ is inverse 
iff for every $m \in M$ there exists one and only one $x \in M$ such that 
$m x m = m$ and $x = x m x$.)

We consider the submonoids ${\it totInv}_{k,1}$ and ${\it surInv}_{k,1}$ 
of ${\it Inv}_{k,1}$, described by bijections 
$P \to Q$ where $P, Q \subset A^*$ are two finite prefix codes with $P$, 
respectively $Q$ maximal. The (unique) inverses of elements in 
${\it totInv}_{k,1}$ are in ${\it surInv}_{k,1}$, and vice versa, so these
submonoids of ${\it Inv}_{k,1}$ are not regular monoids.
We have \ ${\it totInv}_{k,1} \cap {\it surInv}_{k,1} \ = \ G_{k,1}$ \ (the
Thompson-Higman group).

It is easy to see that for all $n>0$, $M_{k,1}$ contains the symmetric monoids
${\it PF}_{k^n}$ of all partial functions on $k^n$ elements, represented by
all elements of $M_{k,1}$ with a table $P \to Q$ where $P, Q \subseteq A^n$.
Hence $M_{k,1}$ contains all finite monoids.                 
Similarly, ${\it tot}M_{k,1}$ contains the symmetric monoids $F_{k^n}$ of 
all total functions on $k^n$ elements.
And ${\it Inv}_{k,1}$ contains ${\mathcal I}_{k^n}$ (the finite symmetric 
inverse monoid of all injective partial functions on $A^n$).


\subsection{Cuntz algebras and Thompson-Higman monoids}
 
All the monoids, inverse monoids, and groups, defined above, are submonoids of
the multiplicative part of the Cuntz algebra ${\cal O}_k$.

The Cuntz algebra ${\cal O}_k$, introduced by Dixmier \cite{Dixmier} (for 
$k=2$) and Cuntz \cite{Cuntz}, is a $k$-generated star-algebra (over the 
field of complex numbers) with identity element {\bf 1} and zero {\bf 0}, 
given by the following finite presentation. 
The generating set is $A = \{a_1, \ldots, a_k\}$.
Since this is defined as a star-algebra, we automatically have the
star-inverses $\{ \overline{a}_1, \ldots, \overline{a}_k\}$; for 
clarity we use overlines rather than stars.

\smallskip

\noindent Relations of the presentation: \

$\overline{a}_i a_i = {\bf 1}$, \ \ for $i = 1, \ldots, k$;

$\overline{a}_i a_j = {\bf 0}$, \ \ when $i \neq j$, $1 \leq i, j \leq k$;

$a_1 \overline{a}_1 + \ldots + a_k \overline{a}_k = {\bf 1}$.

\smallskip

\noindent It is easy to verify that this defines a star-algebra.
The Cuntz algebras are actually C$^*$-algebras with many remarkable 
properties (proved in \cite{Cuntz}), but here we only need them
as star-algebras, without their norm and Cauchy completion.

\medskip

In \cite{BiThomps} and independently in \cite{Nekrash} it was proved that
the Thompson-Higman group $G_{k,1}$ is the subgroup of ${\cal O}_k$ 
consisting of the elements that have an expression of the form 
 \ $\sum_{x \in P} f(x) \ \overline{x}$ \
where we require the following: $P$ and $Q$ range over all finite maximal
prefix codes over the alphabet $\{a_1, \ldots, a_k\}$, and $f$ is any
bijection $P \to Q$.  Another proof is given in \cite{Hughes}.
More generally we also have:

\begin{thm} \label{thomps_monoid_in_cuntz_algebra} \   
The Thompson-Higman monoid $M_{k,1}$ is a submonoid of the multiplicative 
part of the Cuntz algebra ${\cal O}_k$.
\end{thm}
{\bf Proof outline.} 
The Thompson-Higman partial function monoid $M_{k,1}$ is the set of all 
elements of ${\cal O}_k$ that have an expression of the form \  
$\sum_{x \in P} f(x) \ \overline{x}$ \ where
$P \subset A^*$ ranges over all finite prefix codes, and $f$ ranges over 
functions  $P \to A^*$. 

The details of the proof are very similar to the proofs in 
\cite{BiThomps, Nekrash}; the definition of {\it essential} restriction (and
extension) and Proposition \ref{extension_formula} insure that the same proof
goes through. \ \ \ $\Box$

\medskip

The embeddability into the Cuntz algebra is a further justification of the 
definitional choices that we made for the Thompson-Higman monoid $M_{k,1}$.


\section{Structure and simplicity of the Thompson-Higman monoids}

\noindent We give some structural properties of the Thompson-Higman monoids;
in particular, we show that $M_{k,1}$ and ${\it Inv}_{k,1}$ are simple for 
all $k$. 


\subsection{Group of units, $J$-relation, simplicity}

By definition, the group of units of a monoid $M$ is the set of invertible 
elements (i.e., the elements $u \in M$ for which there exists $x \in M$ 
such that $xu = ux = {\bf 1}$, where {\bf 1} is the identity element of $M$). 

\begin{pro} \label{group_of_units} \
The Thompson-Higman group $G_{k,1}$ is the {\em group of units} of the
monoids $M_{k,1}$, ${\it tot}M_{k,1}$, ${\it sur}M_{k,1}$,
${\it totsur}M_{k,1}$, and ${\it Inv}_{k,1}$.
\end{pro}
{\bf Proof.} It is obvious that the groups of units of the above monoids
contain $G_{k,1}$. Conversely, we want to show that that if
$\varphi \in M_{k,1}$ (and in particular, if $\varphi$ is in one of the 
other monoids) and if $\varphi$ has a left inverse and a right inverse,
then $\varphi \in G_{k,1}$.

First, it follows that $\varphi$ is injective, i.e., $\varphi \in $
${\it Inv}_{k,1}$. Indeed, existence of a left inverse implies that for
some $\alpha \in M_{k,1}$ we have $\alpha \ \varphi = {\bf 1}$; hence,
if $\varphi(x_1) = \varphi(x_2)$ then $x_1 = \alpha \ \varphi(x_1) = $
$\alpha \ \varphi(x_2) = x_2$.

Next, we show that ${\rm domC}(\varphi)$ is a {\em maximal} prefix code,
hence $\varphi \in {\it totInv}_{k,1}$. Indeed, we can again consider
$\alpha \in M_{k,1}$ such that $\alpha \ \varphi = {\bf 1}$. For any 
essential restriction of {\bf 1} the domain code is a maximal prefix code,
hence ${\rm domC}(\alpha \circ \varphi)$ is maximal (where $\circ$ denotes
functional composition). Moreover, ${\rm domC}(\alpha \circ \varphi)$ is 
also contained in the domain code of some restriction of $\varphi$, since  
$\varphi(x)$ must be defined when $\alpha \circ \varphi(x)$ is defined. 
Hence domC($\varphi'$), for some restriction $\varphi'$ of $\varphi$, is a 
maximal prefix code; it follows that domC($\varphi$) is a maximal prefix code.

If we apply the reasoning of the previous paragraph to $\varphi^{-1}$ 
(which exists since we saw that $\varphi$ is injective), we conclude
that ${\rm domC}(\varphi^{-1}) = {\rm imC}(\varphi)$ is a maximal prefix code.
Thus, $\varphi \in {\it surInv}_{k,1}$.


We proved that if $\varphi$ has a left inverse and a right inverse then
$\varphi \in {\it totInv}_{k,1} \cap {\it surInv}_{k,1}$.  Since
${\it totInv}_{k,1} \cap {\it surInv}_{k,1} = G_{k,1}$ we conclude that
$\varphi \in G_{k,1}$.
 \ \ \ $\Box$

\bigskip

\noindent We now characterize some of the Green relations of $M_{k,1}$ and 
of ${\it Inv}_{k,1}$, and we prove simplicity.

\medskip

By definition, two elements $x,y$ of a monoid $M$ are {\em $J$-related}
(denoted $x \equiv_J y$) iff $x$ and $y$ belong to exactly the same ideals
of $M$. More generally, the {\em $J$-preorder} of $M$ is defined as follows:
$x \leq_J y$ iff $x$ belongs to every ideal that $y$ belongs to. 
It is easy to see that $x \equiv_J y$ iff $x \leq_J y$ and $y \leq_J x$;
moreover, $x \leq_J y$ iff there exist $\alpha, \beta \in M$ such that 
$x = \alpha y \beta$.
A monoid
$M$ is called {\em $J$-simple} iff $M$ has only one $J$-class (or
equivalently, $M$ has only one ideal, namely $M$ itself). A monoid $M$ is
called {\em $0$-$J$-simple} iff $M$ has exactly two $J$-classes, one of
which consist of just a zero element (equivalently, $M$ has only two ideals,
one of which is a zero element, and the other is $M$ itself).
See \cite{CliffPres, Grillet} for more information on the $J$-relation.
Cuntz \cite{Cuntz} proved that the multiplicative part of the $C^*$-algebra
${\cal O}_k$ is a $0$-$J$-simple monoid, and that as an algebra ${\cal O}_k$
is simple. We will now prove similar results for the Thompson-Higman monoids.

\begin{pro} \label{0Jsimple}
The inverse monoid ${\it Inv}_{k,1}$ and the monoid $M_{k,1}$ are
$0$-$J$-simple.
The monoid ${\it tot}M_{k,1}$ is $J$-simple.
\end{pro}
{\bf Proof.} Let $\varphi \in M_{k,1}$ (or $\in {\it Inv}_{k,1}$). When
$\varphi$ is not the empty map there are $x_0, y_0 \in A^*$ such that
$y_0 = \varphi(x_0)$. Let us define $\alpha, \beta \in {\it Inv}_{k,1}$ by
the tables $\alpha = \{(\varepsilon \mapsto x_0)\}$ and
$\beta = \{(y_0 \mapsto \varepsilon)\}$.
Recall that $\varepsilon$ denotes the empty word.
Then \ $\beta \, \varphi \, \alpha(.) = $
$\{(\varepsilon \mapsto \varepsilon)\} = {\bf 1}$. So, every non-zero element
of $M_{k,1}$ (and of ${\it Inv}_{k,1}$) is in the same $J$-class as the
identity element.

In the case of ${\it tot}M_{k,1}$ we can take
$\alpha = \{(\varepsilon \mapsto x_0)\}$ as before (since the domain code
of $\alpha$ is $\{ \varepsilon\}$, which is a maximal prefix code), and we
take $\beta' : Q \mapsto \{\varepsilon\}$ (i.e., the map that sends
every element of $Q$ to $\varepsilon$), where $Q$ is any finite
maximal prefix code containing $y_0$. Then again, \
$\beta' \, \varphi \, \alpha(.) = $
$\{(\varepsilon \mapsto \varepsilon)\} = {\bf 1}$.  \ \ \ $\Box$

\bigskip

\noindent Thompson proved that $V$ ($= G_{2,1}$) is a simple group; Higman 
proved more generally that when $k$ is even then $G_{k,1}$ is simple, and 
when $k$ is odd then $G_{k,1}$ contains a simple normal subgroup of index 2.
We will show next that in the monoid case we have {\em simplicity for all} 
$k$ (not only when $k$ is even). For a monoid $M$, ``simple'', or more 
precisely, {\em ``congruence-simple''} is defined to mean that the only 
congruences on $M$ are the trivial congruences (i.e., the equality relation, 
and the congruence that lumps all elements of $M$ into one congruence class).

\begin{thm} \label{simple}
The Thompson-Higman monoids ${\it Inv}_{k,1}$ and $M_{k,1}$ are
congruence-simple for all $k$.
\end{thm}  
{\bf Proof.} \ Let $\equiv$ be any congruence on $M_{k,1}$ that is not
the equality relation. We will show that then the whole monoid is
congruent to the empty map {\bf 0}. We will make use of
$0$-$\cal J$-simplicity.

\smallskip

\noindent {\sf Case 0:} \
Assume that $\Phi \equiv $ {\bf 0} for some element $\Phi \neq $
{\bf 0} of $M_{k,1}$. Then for all $\alpha, \beta \in M_{k,1}$ we have
obviously $\alpha \, \Phi \, \beta \equiv $ {\bf 0}. Moreover, by
$0$-$\cal J$-simplicity of $M_{k,1}$ we have \ $M_{k,1}$
$= \{ \alpha \, \Phi \, \beta : \alpha, \beta \in M_{k,1}\}$
 \ since $\Phi \neq 0$. Hence in this case all elements of $M_{k,1}$ are
congruent to {\bf 0}.

\smallskip

\noindent For the remainder we suppose that $\varphi \equiv \psi$ and
$\varphi \neq \psi$, for some elements $\varphi, \psi$ of
$M_{k,1} - \{ {\bf 0} \}$. 

For a right ideal $R \subseteq A^*$ generated by a prefix code $P$
we call $P A^{\omega}$ the set of {\em ends} of $R$.
We call two right ideals $R_1, R_2$ {\em essentially equal}
iff $R_1$ and $R_2$ have the same ends, and we denote this by
$R_1 =_{\sf ess} R_2$. This is equivalent to the following property:
Every right ideal that intersects $R_1$ also intersects $R_2$, and vice 
versa (see \cite{BiBern} and \cite{BiCongr}).

\smallskip

\noindent {\sf Case 1:}
 \ ${\rm Dom}(\varphi) \neq_{\sf ess} {\rm Dom}(\psi)$.

Then there exists $x_0 \in A^*$ such that
$x_0A^* \subseteq {\rm Dom}(\varphi)$, but
 \ ${\rm Dom}(\psi) \cap x_0A^* = \varnothing$;  or, vice versa,
there exists $x_0 \in A^*$ such that
$x_0A^* \subseteq {\rm Dom}(\psi)$, but
 \ ${\rm Dom}(\varphi) \cap x_0A^* = \varnothing$.  Let us assume the former.
Letting $\beta = (x_0 \mapsto x_0)$, we have
$\varphi \, \beta(.) = (x_0 \mapsto \varphi(x_0))$.
We also have $\psi \, \beta(.) = $ {\bf 0}, since
$x_0A^* \cap {\rm Dom}(\psi) = \varnothing$.
So, $\varphi \, \beta \equiv \psi \, \beta = {\bf 0}$, but
$\varphi \, \beta \neq$ {\bf 0}.
Hence case 0, applied to $\Phi = \varphi \, \beta$,
implies that the entire monoid $M_{k,1}$ is congruent to {\bf 0}.

\smallskip

\noindent {\sf Case 2.1:}
 \ ${\rm Im}(\varphi) \neq_{\sf ess} {\rm Im}(\psi)$
 \ and \ ${\rm Dom}(\varphi) =_{\sf ess} {\rm Dom}(\psi)$.

Then there exists $y_0 \in A^*$ such that
$y_0A^* \subseteq {\rm Im}(\varphi)$, but
 \ ${\rm Im}(\psi) \cap y_0A^* = \varnothing$;  or, vice versa,
$y_0A^* \subseteq {\rm Im}(\psi)$, but
 \ ${\rm Im}(\varphi) \cap y_0A^* = \varnothing$.
Let us assume the former. Let $x_0 \in A^*$ be such that
$y_0 = \varphi(x_0)$. Then
 \ $(y_0 \mapsto y_0) \circ \varphi \circ (x_0 \mapsto x_0)$
$ \ = \ (x_0 \mapsto y_0)$.

On the other hand,
 \ $(y_0 \mapsto y_0) \circ \psi \circ (x_0 \mapsto x_0) \ = \ {\bf 0}$.
Indeed, if $x_0A^* \cap {\rm Dom}(\psi) = \varnothing$ then for all
$w \in A^* : $ \ $\psi \circ (x_0 \mapsto x_0)(x_0w) \ = \ \psi(x_0w)$
$ \ = \ $ $\varnothing$.
And if $x_0A^* \cap {\rm Dom}(\psi) \neq  \varnothing$ then
for those $w \in A^*$ such that $x_0w \in {\rm Dom}(\psi)$ we have
 \ $(y_0 \mapsto y_0) \circ \psi \circ (x_0 \mapsto x_0)(x_0w) \ = \ $
$(y_0 \mapsto y_0)(\psi(x_0w)) \ = \ \varnothing$, since
${\rm Im}(\psi) \cap y_0A^* = \varnothing$.  Now case 0 applies to
${\bf 0} \neq \Phi = $
$(y_0 \mapsto y_0) \circ \varphi \circ (x_0 \mapsto x_0)$
$ \equiv {\bf 0}$;
hence all elements of $M_{k,1}$ are congruent to {\bf 0}.

\smallskip

\noindent {\sf Case 2.2:}
  \ ${\rm Im}(\varphi) =_{\sf ess} {\rm Im}(\psi)$
 \ and \ ${\rm Dom}(\varphi) =_{\sf ess} {\rm Dom}(\psi)$.

Then after restricting $\varphi$ and $\psi$ to 
${\rm Dom}(\varphi) \cap {\rm Dom}(\psi)$ ($ =_{\sf ess}$
${\rm Dom}(\varphi) =_{\sf ess} {\rm Dom}(\psi)$), we have:  
$\, {\rm domC}(\varphi) = {\rm domC}(\psi)$, and
 there exist $x_0 \in {\rm domC}(\varphi) = {\rm domC}(\psi)$ and
$y_0 \in {\rm Im}(\varphi)$, $y_1 \in {\rm Im}(\psi)$ such
that $\varphi(x_0) = y_0 \neq y_1 = \psi(x_0)$. We have two sub-cases.

\smallskip

\noindent {\sf Case 2.2.1:} \ $y_0$ and $y_1$ are not prefix-comparable.

Then
 \ $(y_0 \mapsto y_0) \circ \varphi \circ (x_0 \mapsto x_0)$
$ \ = \ (x_0 \mapsto y_0)$.

On the other hand,
 \ $(y_0 \mapsto y_0) \circ \psi \circ (x_0 \mapsto x_0)(x_0 w) \ = $
 \ $(y_0 \mapsto y_0)(y_1 w) \ = \ \varnothing$ \ for all $w \in A^*$ \
(since $y_0$ and $y_1$ are not prefix-comparable).  So
 \ $(y_0 \mapsto y_0) \circ \psi \circ (x_0 \mapsto x_0) \ = \ \ {\bf 0}$.
Hence case 0 applies to \ ${\bf 0} \neq\Phi \ = \ $
$(y_0 \mapsto y_0) \circ \varphi \circ (x_0 \mapsto x_0)$
$\equiv {\bf 0}$.

\smallskip

\noindent {\sf Case 2.2.2:} \ $y_0$ is a prefix of $y_1$, and
$y_0 \neq y_1$. (The case where $y_0$ is a prefix of $y_1$ is similar.)

Then $y_1 = y_0 a u_1$ for some $a \in A$, $u_1 \in A^*$.
Letting $b \in A - \{a\}$, and $y_2 = y_0 b$, we obtain a string $y_2$ that
is not prefix-comparable with $y_1$.
Now, $\, (y_2 \mapsto y_2) \circ \varphi \circ (x_0 \mapsto x_0)(x_0v_2)$
$ \ = \ $ $(y_2 \mapsto y_2)(y_0 v_2) \ = \ y_2$.
But for all $w \in A^*$, 
$ \, (y_2 \mapsto y_2) \circ \psi \circ (x_0 \mapsto x_0)(x_0w) \ = \ $
 $(y_2 \mapsto y_2)(y_1 w) \ = \ \varnothing$, since
$y_2$ and $y_1$ are not prefix-comparable.
Thus, case 0 applies to \ ${\bf 0} \neq \Phi = $
$(y_2 \mapsto y_2) \circ \varphi \circ (x_0 \mapsto x_0)$
$\equiv {\bf 0}$.

\smallskip

The same proof works for ${\it Inv}_{k,1}$ since all the multipliers used
in the proof (of the form $(u \mapsto v)$ for some $u,v \in A^*$) belong to
${\it Inv}_{k,1}$.
 \ \ \ $\Box$

\bigskip


\subsection{$D$-relation}

Besides the $J$-relation and the $J$-preorder, based on ideals, there are 
the $R$- and $L-$relations and $R$- and $L-$preorders, based on right (or 
left) ideals.
Two elements $x,y \in M$ are {\em $R$-related} (denoted $x \equiv_R y$) iff
$x$ and $y$ belong to exactly the same right ideals of $M$. The
{\em $R$-preorder} is defined as follows: $x \leq_R y$ iff $x$ belongs to
every right ideal that $y$ belongs to. 
It is easy to see that $x \equiv_R$ iff $x \leq_R y$ and $y \leq_R x$; also,
$x \leq_R y$ iff there exists $\alpha \in M$ such that $x = y \alpha$. 
In a similar way one defines $\equiv_L$ and $\leq_L$.
Finally, there is the {\em $D$-relation} of $M$, which is defined as follows:
$x \equiv_D y$ iff there exists $s \in M$ such that $x \equiv_R s \equiv_L y$;
this is easily seen to be equivalent to saying that there exists $t \in M$
such that $x \equiv_L t \equiv_R y$. For more information on these definitions
see for example \cite{CliffPres, Grillet}.

The $D$-relation of $M_{k,1}$ and ${\it Inv}_{k,1}$ has an interesting
characterization, as we shall prove next.
We will represent all elements of $M_{k,1}$ by tables of the from
$\varphi: P \to Q$, where both $P$ and $Q$ are finite prefix codes over $A$
(with $|A| = k$). For such a table we also write $P = {\rm domC}(\varphi)$
(the domain code of $\varphi$) and $Q = {\rm imC}(\varphi)$ (the image code
of $\varphi$). In general, tables of elements of $M_{k,1}$ have the form
$P \to S$, where $P$ is a finite prefix code and $S$ is a finite set; but
by using essential restrictions, if necessary, every element of $M_{k,1}$ can
be given a table $P \to Q$, where both $P$ and $Q$ are finite prefix codes.

Note the following invariants with respect to essential restrictions:

\begin{pro} \label{essent_ext_restr_mod_k_1} \
Let $\varphi_1: P_1 \to Q_1$ be a table for an element of $M_{k,1}$, where
$P_1,Q_1 \subset A^*$ are finite prefix codes.  Let $\varphi_2: P_2 \to Q_2$ 
be another finite table for the {\em same} element of $M_{k,1}$,
obtained from the table $\varphi_1$ by an essential restriction.
Then $P_2, Q_2 \subset A^*$ are finite prefix codes and we have

\smallskip

 \ \ \ \ \ $|P_1| \equiv |P_2|$ \ {\rm mod} $(k-1)$ \ \ \ and

\smallskip

 \ \ \ \ \ $|Q_1| \equiv |Q_2|$ \ {\rm mod} $(k-1)$.

\smallskip

\noindent These modular congruences also hold for essential extensions,
provided that we only extend to tables in which the image is a prefix 
code.
\end{pro}
{\bf Proof.} An essential restriction consists of a finite sequence of
essential restriction steps; an essential restriction step consists of 
replacing a table entry $(x,y)$ of $\varphi_1$ by
$\{(xa_1, ya_1), \ldots, (xa_k, ya_k)\}$ (according to Proposition
\ref{extension_formula}).
For a finite prefix code $Q \subset A^*$, and $q \in Q$, the finite set
$(Q - \{q\}) \cup \{qa_1, \ldots, qa_k\}$ is also a prefix code, as is 
easy to prove.
In this process, the cardinalities change as follows: \ $|P_1|$ becomes
$|P_1| - 1 + k$ and $|Q_1|$ becomes $|Q_1| - 1 + k$. Indeed (looking at
$Q_1$ for example), first an element $y$ is removed from $Q_1$, then the
$k$ elements $\{ya_1, \ldots, ya_k\}$ are added. The elements $ya_i$ that
are added are all different from the elements that are already present in
$Q_1 - \{y\}$; in fact, more strongly, $ya_i$ and the elements of
$Q_1 - \{y\}$ are not prefixes of each other. 
 \ \ \ $\Box$

\medskip

As a consequence of Prop.\ \ref{essent_ext_restr_mod_k_1} it makes sense,
for any $\varphi \in M_{k,1}$, to talk about $|{\rm domC}(\varphi)|$ and
$|{\rm imC}(\varphi)|$ as elements of ${\mathbb Z}_{k-1}$, independently of
the representation of $\varphi$ by a right-ideal homomorphism.

\begin{thm} \label{D_relation_characteris}
For any non-zero elements $\varphi, \psi$ of $M_{k,1}$ (or of
${\it Inv}_{k,1}$) the $D$-relation is characterized as follows:

\smallskip

 \ \ \ $\varphi \equiv_D \psi$ \ \ \ iff \ \ \  
$|{\rm imC}(\varphi)| \equiv |{\rm imC}(\psi)|$ \ {\rm mod} $(k-1)$.

\smallskip

\noindent
Hence, $M_{k,1}$ and ${\it Inv}_{k,1}$ have $k-1$ non-zero $D$-classes.
In particular, $M_{2,1}$ and ${\it Inv}_{2,1}$ are $0$-$D$-simple (also
called {\em $0$-bisimple}).
\end{thm}
The proof of Theorem \ref{D_relation_characteris} uses several Lemmas.

\begin{lem} \label{card_max_pref_code}
{\rm (\cite{BiCoNP} Lemma 6.1; Arxiv version of \cite{BiCoNP} Lemma 9.9).}  \   
For every finite alphabet $A$ and every integer $i \geq 0$ there exists a 
maximal prefix code of cardinality \ $1 + (|A|-1) \, i$.
And every finite maximal prefix code over $A$ has cardinality 
 \ $1 + (|A|-1) \, i$, for some integer $i \geq 0$.

It follows that when $|A| = 2$, there are finite prefix codes over $A$ of 
every finite cardinality. \ \ \ \ \ $\Box$
\end{lem}
As a consequence of this Lemma we have for all $\varphi \in G_{k,1}$: \
$\|\varphi\| \equiv 1$ mod $(k-1)$. Thus, except for the Thompson group
$V$ (when $k=2$), there is a constraint on the table size of the elements
of the group.

In the following ${\rm id}_Q$ denotes the element of
${\it Inv}_{k,1}$ given by the table $\{(x \mapsto x) : x \in Q\}$  where
$Q \subset A^*$ is any finite prefix code.

\begin{lem} \label{D_if_imC_same_card} \
{\bf (1)}
For any $\varphi \in M_{k,1}$ (or $\in {\it Inv}_{k,1}$) with table
$P \to Q$ (where $P, Q$ are finite prefix codes) we have: \
$\varphi \equiv_R {\rm id}_Q$.  \\
{\bf (2)}
If $S, T$ are finite prefix codes with $|S| = |T|$ then \
${\rm id}_S \equiv_D {\rm id}_T$. \\
{\bf (3)}   If $\varphi_1: P_1 \to Q_1$ and $\varphi_2: P_ 2 \to Q_2$ are
such that $|Q_1| = |Q_2|$ then \ $\varphi_1 \equiv_D \varphi_2$.
\end{lem}
{\bf Proof.}
{\bf (1)} Let $P' \subseteq P$ be a set of representatives modulo $\varphi$
(i.e., we form $P'$ by choosing one element in every set
$\varphi^{-1} \varphi(x)$ as $x$ ranges over $P$). So, $|P'| = |Q|$.
Let $\alpha \in {\it Inv}_{k,1}$ be given by a table $Q \to P'$; the exact
map does not matter, as long as $\alpha$ is bijective.
Then $\varphi \circ \alpha(.)$ is a permutation of $Q$, and
$\varphi \circ \alpha \equiv_R $
$\varphi \circ \alpha \circ (\varphi \circ \alpha)^{-1} = {\rm id}_Q$.

Now, $\varphi \geq_R \varphi \circ \alpha \geq_R $
$\varphi \circ \alpha \circ (\varphi \circ \alpha)^{-1} \circ \varphi = $
${\rm id}_Q \circ \varphi = \varphi$, hence 
$\varphi \equiv_R \varphi \circ \alpha$ \ ($\equiv_R {\rm id}_Q$).

\smallskip

\noindent {\bf (2)} Let $\alpha: S \to T$ be a bijection (which exists
since $|S| = |T|$); so $\alpha$ represents an element of ${\it Inv}_{k,1}$.
Then $\alpha = \alpha \circ {\rm id}_S(.)$ and
${\rm id}_S = \alpha^{-1} \circ \alpha(.)$; hence,
$\alpha \equiv_L {\rm id}_S$.

Also, $\alpha = {\rm id}_T \circ \alpha(.)$ and
${\rm id}_T = \alpha \circ \alpha^{-1}(.)$; hence,
$\alpha \equiv_R {\rm id}_T$.
Thus, ${\rm id}_S \equiv_L \alpha \equiv_R {\rm id}_T$.

\smallskip

\noindent {\bf (3)} If $|Q_1| = |Q_2|$ then
${\rm id}_{Q_1} \equiv_D {\rm id}_{Q_2}$ by (2). Moreover,
$\varphi_1 \equiv_D {\rm id}_{Q_1}$ and $\varphi_2 \equiv_D {\rm id}_{Q_2}$
by (1). The result follows by transitivity of $\equiv_D$.
 \ \ \ $\Box$

\begin{lem} \label{special_pref_code} \
{\bf (1)}
For any $m \geq k$ let $i$ be the residue of $m$ modulo $k-1$ in the range
$2 \leq i \leq k$, and let us write $m = i + (k-1)j$, for some $j \geq 0$.
Then there exists a prefix code $Q_{i,j}$ of cardinality $|Q_{i,j}| = m$,
such that ${\rm id}_{Q_{i,j}}$ is an essential restriction of
${\rm id}_{\{a_1, \ldots,a_i\}}$. Hence, \
 ${\rm id}_{Q_{i,j}} = {\rm id}_{\{a_1, \ldots,a_i\}}$ \ as elements of
${\it Inv}_{k,1}$. \\  
{\bf (2)} In $M_{k,1}$ and in ${\it Inv}_{k,1}$ we have \
${\rm id}_{\{a_1\}} \equiv_D {\rm id}_{\{a_1, \ldots, a_k\}} = {\bf 1}$.
\end{lem}
{\bf Proof.} {\bf (1)} For any $m \geq k$ there exist $i, j \geq 0$ such that
$1 \leq i \leq k$ and $m = i + (k-1)j$.  
We consider the prefix code

\smallskip

$Q_{i,j} \ = \ \{a_2, \ldots,a_i\} \ \ \cup \ \ $
$\bigcup_{r=1}^{j-1} a_1^r (A-\{a_1\}) \ \ \cup \ \ a_1^j A$.

\smallskip

\noindent It is easy to see that $Q_{i,j}$ is a prefix code, which is maximal
iff $i = k$; see Fig.\ 1 below. Clearly, $|Q_{i,j}| = i + (k-1)j$.
Since $Q_{i,j}$ contains $a_1^j A$, we can perform an essential
extension of ${\rm id}_{Q_{i,j}}$ by replacing the table entries \
$\{(a_1^j a_1, a_1^j a_1), (a_1^j a_2, a_1^j a_2), \ldots, $
$(a_1^j a_k,a_1^j a_k)\}$ \ by $(a_1^j, a_1^j)$. This replaces $Q_{i,j}$
by $Q_{i,j-1}$. So, ${\rm id}_{Q_{i,j}}$ can be essentially extended to
${\rm id}_{Q_{i,j-1}}$.  By repeating this we find that ${\rm id}_{Q_{i,j}}$
is the same element (in $M_{k,1}$ and in ${\it Inv}_{k,1}$) as
${\rm id}_{Q_{i,0}} = {\rm id}_{\{a_1, \ldots, a_i\}}$.

\smallskip

\noindent {\bf (2)} By essential restriction, \
${\rm id}_{\{a_1\}} = {\rm id}_{\{a_1a_1, a_1a_2, \ldots, a_1a_k\}}$, in
$M_{k,1}$ and in ${\it Inv}_{k,1}$. And by Lemma
\ref{D_if_imC_same_card}(2), \
${\rm id}_{\{a_1a_1, a_1a_2, \ldots, a_1a_k\}} \equiv_D $
${\rm id}_{\{a_1, \ldots, a_k\}}$; the latter, by essential extension, is
{\bf 1}.
 \ \ \ $\Box$

\bigskip

\unitlength=0.90mm \special{em:linewidth 0.4pt}
\linethickness{0.4pt}

\begin{picture}(150,100)

\put(120,90){\makebox(0,0)[cc]{$\varepsilon$}}
\put(123,87){\line(1,-2){2}}
\put(136,79){\makebox(0,0)[cc]{$\{a_2,\ldots,a_i\}$}}
\put(124,87){\line(2,-1){8}}
\put(128,83){\makebox(0,0)[cc]{$\ldots$}}

\put(117,87){\line(-1,-2){2.1}}
\put(114,80){\makebox(0,0)[cc]{$a_1$}}

\put(117,77){\line(1,-2){2.2}}
\put(133,70){\makebox(0,0)[cc]{$a_1 \, (A-\{a_1\})$}}
\put(118,77){\line(2,-1){8}}
\put(122,73){\makebox(0,0)[cc]{$\ldots$}}

\put(111,77){\line(-1,-2){2}}
\put(107,70){\makebox(0,0)[cc]{.}}
\put(106,68){\makebox(0,0)[cc]{.}}
\put(105,66){\makebox(0,0)[cc]{.}}
\put(103,62){\line(-1,-2){2}}

\put(99,54){\makebox(0,0)[cc]{$a_1^r$}}
\put(102,50){\line(1,-2){2}}
\put(119,43){\makebox(0,0)[cc]{$a_1^r \, (A-\{a_1\})$}}
\put(103,50){\line(2,-1){8}}
\put(107,46){\makebox(0,0)[cc]{$\ldots$}}

\put(96,50){\line(-1,-2){2}}
\put(92,43){\makebox(0,0)[cc]{.}}
\put(91,41){\makebox(0,0)[cc]{.}}
\put(90,39){\makebox(0,0)[cc]{.}}
\put(88,35){\line(-1,-2){2}}

\put(84,26){\makebox(0,0)[cc]{$a_1^{j-1}$}}
\put(85,22){\line(1,-2){2.7}}
\put(104,13){\makebox(0,0)[cc]{$a_1^{j-1} \, (A-\{a_1\})$}}
\put(86,22){\line(2,-1){10}}
\put(91,17){\makebox(0,0)[cc]{$\ldots$}}

\put(79,22){\line(-1,-2){2.2}}
\put(73,13){\makebox(0,0)[cc]{$a_1^j$}}

\put(71,9){\line(-1,-2){2.2}}
\put(75,9){\line(1,-2){2.2}}
\put(73,5){\makebox(0,0)[cc]{$\ldots$}}

\put(73,1){\makebox(0,0)[cc]{$A$}}

\put(20,-6){\makebox(0,0)[cc]{{\sf Fig.\ 1: \ The prefix tree of $Q_{i,j}$.}}}

\end{picture}

\bigskip

\bigskip

\begin{lem} \label{Green_rel_in_Inv_from_M} \   
For all $\varphi, \psi \in {\it Inv}_{k,1}$: If \  
$\varphi \geq_{L(M_{k,1})} \psi$, where $\geq_{L(M_{k,1})}$ is the 
$L$-preorder of  $M_{k,1}$, then \ $\varphi \geq_{L(I_{k,1})} \psi$, where 
$\geq_{L(I_{k,1})}$ is the $L$-preorder of ${\it Inv}_{k,1}$.

The same holds with $\geq_L$ replaced by $\equiv_L$, $\geq_R$, $\equiv_R$,
$\equiv_D$, $\geq_J$ and $\equiv_J$. 
\end{lem}
{\bf Proof.} If $\psi = \alpha \, \varphi$ for some $\alpha \in M_{k,1}$
then let us define $\alpha'$ by \ 
$\alpha' = \alpha \ {\rm id}_{{\rm Im}(\varphi)}$. Then we have:
 \ $\psi \ \varphi^{-1} = \alpha \ \varphi \ \varphi^{-1} = $
$\alpha \ {\rm id}_{{\rm Im}(\varphi)} = \alpha'$, hence 
$\alpha' \in {\it Inv}_{k,1}$ (since $\varphi, \psi \in {\it Inv}_{k,1}$).
Moreover, $\alpha' \, \varphi = $
$\alpha \ {\rm id}_{{\rm Im}(\varphi)} \  \varphi = \alpha \, \varphi = $
$\psi$.
 \ \ \ $\Box$

\bigskip

So far our Lemmas imply that in $M_{k,1}$ and in ${\it Inv}_{k,1}$, every 
non-zero element is $\equiv_D$ to one of the $k-1$ elements 
${\rm id}_{\{a_1, \ldots, a_i\}}$, for $i=1, \ldots, k-1$. 
Moreover the Lemmas show that if two elements of $M_{k,1}$ (or of
${\it Inv}_{k,1}$) are given by tables $\varphi_1: P_1 \to Q_1$ and 
$\varphi_2: P_2 \to Q_2$, where $P_1$, $Q_1$, $P_2$ and $Q_2$ are finite 
prefix codes, then we have: \ If \ $|Q_1| \equiv |Q_2|$ mod $(k-1)$ \ then 
 \ $\varphi_1 \equiv_D \varphi_2$.

We still need to prove the converse of this. It is sufficient to prove the
converse for ${\it Inv}_{k,1}$, by Lemma \ref{Green_rel_in_Inv_from_M} and
because every element of $M_{k,1}$ is $\equiv_D$ to an element of
${\it Inv}_{k,1}$ (namely ${\rm id}_{\{a_1,\ldots,a_i\}}$).

\begin{lem} \label{L_implies_mod_k1} \
Let $\varphi, \psi \in {\it Inv}_{k,1}$.
If $\varphi \equiv_D \psi$ in ${\it Inv}_{k,1}$, then 
 \ $\|\varphi\| \equiv \|\psi\|$ {\rm mod} $(k-1)$.
\end{lem}
{\bf Proof.}
 (1) We first prove that if $\varphi \equiv_L \psi$ then
$|{\rm domC}(\varphi)| \equiv |{\rm domC}(\psi)|$ {\rm mod} $(k-1)$.

By definition, $\varphi \equiv_L \psi$ iff $\varphi = \beta \, \psi$ and 
$\psi = \alpha \, \varphi$ for some $\alpha, \beta \in {\it Inv}_{k,1}$. 
By Lemma \ref{dom_equal_im} \ there are restrictions $\beta'$ and $\psi'$ 
of $\beta$, respectively $\psi$, and an essential restriction $\Phi$ of 
$\varphi$ such that:

\smallskip

 $\Phi = \beta' \circ \psi'$, \ and \
 ${\rm Dom}(\beta') = {\rm Im}(\psi')$.

\smallskip

\noindent It follows that ${\rm Dom}(\Phi) \subseteq {\rm Dom}(\psi')$,
since if $\psi'(x)$ is not defined then $\Phi(x) = \beta' \circ \psi'(x)$
is not defined either.
Similarly, there is an essential restriction $\Psi$ of $\psi$ and a
restriction $\varphi'$ of $\varphi$ and such that
${\rm Dom}(\Psi) \subseteq {\rm Dom}(\varphi')$.

Thus, the restriction of both $\varphi$ and $\psi$ to the intersection
${\rm Dom}(\Phi) \cap {\rm Dom}(\Psi)$ yields restrictions $\varphi''$, 
respectively $\psi''$ such that ${\rm Dom}(\varphi'') = {\rm Dom}(\psi'')$.

\smallskip

\noindent {\sf Claim:} 
 \ $\varphi''$ and $\psi''$ are essential restrictions of $\varphi$,
respectively $\psi$. 

Indeed, every
right ideal $R$ of $A^*$ that intersects ${\rm Dom}(\psi)$ also intersects
${\rm Dom}(\Psi)$ (since $\Psi$ is an essential restriction of $\psi$). 
Since ${\rm Dom}(\Psi) \subseteq {\rm Dom}(\varphi') \subseteq $
${\rm Dom}(\varphi)$, it follows that $R$ also intersects ${\rm Dom}(\varphi)$.
Moreover, since $\Phi$ is an essential restriction of $\varphi$, $R$ also
intersects ${\rm Dom}(\Phi)$. 
Thus, ${\rm Dom}(\Phi)$ is essential in ${\rm Dom}(\psi)$. 
Since ${\rm Dom}(\Psi)$ is also essential in ${\rm Dom}(\psi)$, it follows
that ${\rm Dom}(\Phi) \cap {\rm Dom}(\Psi)$ is essential in 
${\rm Dom}(\psi)$; indeed, in general, the intersection of two right ideals 
$R_1, R_2$ that are essential in a right ideal $R_3$, is essential in $R_3$
(this is a special case of Lemma \ref{inters_essent_ideals}).    
This means that $\psi''$ is an essential restriction of $\psi$.
Similarly, one proves that $\varphi''$ is an essential restriction of
$\varphi$. \ \ \ [This proves the Claim.]

\smallskip

So, $\varphi''$ and $\psi''$ are essential restrictions such that 
${\rm Dom}(\varphi'') = {\rm Dom}(\psi'')$.
Hence, ${\rm domC}(\varphi'') = {\rm domC}(\psi'')$; Proposition
\ref{essent_ext_restr_mod_k_1} then implies that \
$|{\rm domC}(\varphi)| \equiv |{\rm domC}(\varphi'')|  = $
$|{\rm domC}(\psi'')| \equiv |{\rm domC}(\psi)|$ \ mod \ $(k-1)$.

\smallskip

\noindent (2) Next, let us prove that if $\varphi \equiv_R \psi$ then
$|{\rm imC}(\varphi)| \equiv |{\rm imC}(\psi)|$ {\rm mod} $(k-1)$.
In ${\it Inv}_{k,1}$ we have $\varphi \equiv_R \psi$ iff 
$\varphi^{-1} \equiv_L \psi^{-1}$. Also,
${\rm imC}(\varphi) = {\rm domC}(\varphi^{-1})$. Hence, (2) follows 
from (1). 

\smallskip

\noindent The Lemma now follows from (1) and (2), since for elements of
${\it Inv}_{k,1}$, $|{\rm imC}(\varphi)| = |{\rm domC}(\varphi)| = $
$\|\varphi\|$, and since the $D$-relation is the composite of the $L$-relation
and the $R$-relation.
 \ \ \ $\Box$

\medskip

\noindent {\bf Proof of Theorem \ref{D_relation_characteris}.}
We saw already (in the observations before Lemma \ref{L_implies_mod_k1} and
in the preceding Lemmas) that for $\varphi_1: P_1 \to Q_1$ and
$\varphi_2: P_2 \to Q_2$ (where $P_1$, $Q_1$, $P_2$ and $Q_2$ are non-empty
finite prefix codes) we have: \ If \
$|Q_1| \equiv |Q_2|$ mod $(k-1)$ \ then \ $\varphi_1 \equiv_D \varphi_2$.
In particular, when $|Q_1| \equiv i$ mod $(k-1)$ then
$\varphi_1 \equiv_D {\rm id}_{\{a_1, \ldots, a_i\}}$.

It follows from Lemma \ref{L_implies_mod_k1} that the elements
${\rm id}_{\{a_1, \ldots, a_i\}}$ (for $i=1, \ldots, k-1$) are all in
different $D$-classes.
 \ \ \ $\Box$

\bigskip

So far we have characterized the $D$- and $J$-relations of $M_{k,1}$ and
${\it Inv}_{k,1}$. We leave the general study of the Green relations of 
$M_{k,1}$, ${\it Inv}_{k,1}$, and the other Thompson-Higman monoids for
future work. The main result of this paper, to be proved next, is that
the Thompson-Higman monoids $M_{k,1}$ and ${\it Inv}_{k,1}$ are finitely 
generated and that their word problem over any finite generating set is in 
{\sf P}.


\section{Finite generating sets}

We will show that ${\it Inv}_{k,1}$ and $M_{k,1}$ are finitely generated.
An application of the latter fact is that a finite generating set of 
$M_{k,1}$ can be used to build combinational circuits for finite boolean 
functions that do not have fixed-length inputs or outputs.
In engineering, non-fixed length inputs or outputs make sense, for example,
if the inputs or outputs are handled sequentially, and if the possible input 
strings form a prefix code.

First we need some more definitions about prefix codes. The {\em prefix tree}
of a prefix code $P \subset A^*$ is, by definition, a tree whose vertex set
is the set of all the prefixes of the elements of $P$, and whose edge set is
$\{ (x,xa) : a \in A, \ xa$ is a prefix of some element of $P\}$. The tree is
rooted, with root $\varepsilon$ (the empty word). Thus, the prefix tree of $P$
is a subtree of the tree of $A^*$. The set of leaves of the prefix
tree of $P$ is $P$ itself. The vertices that are not leaves are called 
{\em internal vertices}. We will say more briefly an ``internal vertex
of $P$'' instead of internal vertex of the prefix tree of $P$.
An internal vertex has between 1 and $k$ children; an internal vertex is
called {\em saturated} iff it has $k$ children. 

One can prove easily that a prefix code $P$ is maximal iff every internal 
vertex of the prefix tree of $P$ is saturated. Hence, every prefix code $P$ 
can be embedded in a maximal prefix code (which is finite when $P$ is finite),
obtained by saturating the prefix tree of $P$. Moreover we have:

\begin{lem} \label{ext_to_max_pref_code} \
For any two finite non-maximal prefix codes $P_1, P_2 \subset A^*$
there are finite maximal prefix codes $P'_1, P'_2 \subset A^*$ such that
$P_1 \subset P'_1$, $P_2 \subset P'_2$, and $|P'_1| = |P'_2|$.
\end{lem}
{\bf Proof.} First we saturate $P_1$ and $P_2$ to obtain two maximal prefix
codes $P''_1$ and $P''_2$ such that $P_1 \subset P''_1$, and
$P_2 \subset P''_2$. If $|P''_1| \neq |P''_2|$ (e.g., if $|P''_1| < |P''_2|$)
then $|P''_1|$ and $|P''_2|$ differ by a multiple of $k-1$ (by Prop.\  
\ref{essent_ext_restr_mod_k_1}).
So, in order to make $|P''_1|$ equal to $|P''_2|$ we repeat the following
(until $|P''_1| = |P''_2|$): consider a leaf of the prefix tree of $P''_1$
that does not belong to $P_1$, and attach $k$ children at that leaf; now
this leaf is no longer a leaf, and the net increase in the number of leaves
is $k-1$.
 \ \ \ $\Box$

\begin{lem} \label{transitive_action_G_k1_on_prefixCodes}
Let $P$ and $Q$ be finite prefix codes of $A^*$ with $|P| = |Q|$.
If $P$ and $Q$ are both maximal prefix codes, or if both are non-maximal,
then there is an element of $G_{k,1}$ that maps $P$ onto $Q$.
On the other hand, if one of $P$ and $Q$ is maximal and the other one is not
maximal, then there is no element of $G_{k,1}$ that maps $P$ onto $Q$.
\end{lem}
{\bf Proof.} When $P$ and $Q$ are both maximal then any one-to-one
correspondence between $P$ and $Q$ is an element of $G_{k,1}$.

When $P$ and $Q$ are both non-maximal, we use Lemma \ref{ext_to_max_pref_code}
above to find two maximal prefix codes $P'$ and $Q'$ such that
$P \subset P'$, $Q \subset Q'$, and $|P'| = |Q'|$. Consider now any bijection
from $P'$ onto $Q'$ that is also a bijection from $P$ onto $Q$. This is an
element of $G_{k,1}$.

When $P$ is maximal and $Q$ is non-maximal, then every element
$\varphi \in M_{k,1}$ that maps $P$ onto $Q$ will satisfy
${\rm domC}(\varphi) = P$; since $\varphi$ is onto $Q$, we have
${\rm imC}(\varphi) = Q$. Hence, $\varphi \not\in G_{k,1}$ since
${\rm imC}(\varphi)$ is a non-maximal prefix code. A similar reasoning shows
that no element of $G_{k,1}$ maps $P$ onto $Q$ if $P$ is non-maximal and $Q$
is maximal.
 \ \ \ $\Box$

\medskip

\noindent Notation: For $u,v \in A^*$, the element of ${\it Inv}_{k,1}$
with one-element domain code $\{u\}$ and one-element image code $\{v\}$ is
denoted by $(u \mapsto v)$. When $(u \mapsto v)$ is composed with itself 
$j$ times the resulting element of ${\it Inv}_{k,1}$ is denoted by
$(u \mapsto v)^j$.
\begin{lem} \label{special_generation} \
{\bf (1)} \  For all $j > 0$: \ \
$(a_1 \mapsto a_1a_1)^j \ = \ (a_1 \mapsto a_1^{j+1})$. \\
{\bf (2)} \
Let \ $S = {\{a_1^j a_1, a_1^j a_2, \ldots, a_1^j a_i\}}$, for some
$1 \leq i \leq k-1$, $0 \leq j$. Then ${\rm id}_S$
is generated by the $k+1$ elements
 \ $\{(a_1 \mapsto a_1a_1), \ (a_1a_1 \mapsto a_1)\} \ \cup \ $
$\{{\rm id}_{\{a_1a_1, \ a_1a_2, \ \ldots, \ a_1a_i\}} : 1 \leq i \leq k-1\}$.
   \\
{\bf (3)} \ For all $j \geq 2$: \ \   
$(\varepsilon \mapsto a_1^j)(.) \ = \ (a_1 \mapsto a_1a_1)^{j-1} \, \cdot $
   $(\varepsilon \mapsto a_1)(.)$.
\end{lem}
{\bf Proof.} {\bf (1)} We prove by induction that \
$(a_1 \mapsto a_1a_1)^j = (a_1 \mapsto a_1 a_1^j)$ \ for all $j \geq 1$. \\  
Indeed, \ $(a_1 \mapsto a_1a_1)^{j+1}(.) = (a_1 \mapsto a_1a_1) \, \cdot$
  $(a_1 \mapsto a_1a_1^j)(.)$, \ and by essential restriction this is

\medskip

$\left( \hspace{-.08in} \begin{array}{l|l}
a_1 a_1^j     & \ a_1 w \ \ \ \ (w \in A^j - \{a_1^j\})   \\
a_1a_1 a_1^j  & \ a_1a_1 w
\end{array} \hspace{-.08in} \right)
 \cdot (a_1 \mapsto a_1 a_1^j)(.) $
$ \ = \ (a_1 \mapsto a_1 a_1 a_1^j)(.)$.

\medskip

\noindent {\bf (2)} For
$S = {\{a_1^j a_1, a_1^j a_2, \ldots, a_1^j a_i\}}$ we have

\medskip

${\rm id}_S \ = \ $
$\left( \hspace{-.08in} \begin{array}{l|l|l|l}
a_1 a_1 & a_1a_2 & \ \ldots \ & a_1a_i  \\
a_1^ja_1 & a_1^ja_2 & \ \ldots \ & a_1^ja_i
\end{array} \hspace{-.08in} \right) \cdot $
$\left( \hspace{-.08in} \begin{array}{l|l|l|l}
a_1^j a_1 & a_1^ja_2 & \ \ldots \ & a_1^j a_i  \\
a_1a_1 & a_1a_2 & \ \ldots \ & a_1a_i
\end{array} \hspace{-.08in} \right)(.)$

\smallskip

\noindent      and

\smallskip

$\left( \hspace{-.08in} \begin{array}{l|l|l|l}
a_1 a_1 & a_1a_2 & \ \ldots \ & a_1a_i  \\
a_1^ja_1 & a_1^ja_2 & \ \ldots \ & a_1^ja_i
\end{array} \hspace{-.08in} \right) $

$ \ = \ \left( \hspace{-.08in} \begin{array}{l|l|l|l|l|l|l}
a_1 a_1 & a_1a_2 & \ \ldots \ & a_1a_i & a_1a_{i+1} & \ \ldots \ & a_1a_k \\
a_1^ja_1 & a_1^ja_2 & \ldots & a_1^ja_i & a_1^ja_{i+1} & \ \ldots \ & a_1^j a_k
\end{array}  \right) \cdot $
${\rm id}_{\{a_1a_1,\ a_1a_2,\ \ldots, \ a_1a_i\}}(.)$

\medskip

$ \ = \ (a_1 \mapsto a_1^j) \cdot $
${\rm id}_{\{a_1a_1, \ a_1a_2, \ \ldots, \ a_1a_i\}}$
$ \ = \ $
$(a_1 \mapsto a_1a_1)^{j-1} \cdot $
${\rm id}_{\{a_1a_1, \ a_1a_2, \ \ldots, \ a_1a_i\}}$.

\medskip

\noindent The map id$_{\{a_1a_1\}}$ is redundant as a generator since
$(a_1 a_1 \mapsto a_1 a_1) = (a_1a_1 \mapsto a_1) \, $
$(a_1 \mapsto a_1a_1)(.)$.

\medskip

\noindent {\bf (3)}  By (1) we have \  
$(\varepsilon \mapsto a_1^j) = (a_1 \mapsto a_1^j) \, \cdot$
$ (\varepsilon \mapsto a_1)(.)$, and \ $(a_1 \mapsto a_1^j) = $
$(a_1 \mapsto a_1 a_1)^{j-1}$.
 \ \ \ $\Box$

\begin{thm} \label{Inv_fin_gen} \ \
The inverse monoid ${\it Inv}_{k,1}$ is finitely generated.
\end{thm}
{\bf Proof.} Our strategy for finding a finite generating set for
${\it Inv}_{k,1}$ is as follows: We will use the fact that the
Thompson-Higman group $G_{k,1}$ is finitely generated. Hence, if
$\varphi \in {\it Inv}_{k,1}$, \, $g_1, g_2 \in G_{k,1}$, and if
$g_2 \varphi g_1$ can be expressed as a product $p$ over a fixed finite set
of elements of ${\it Inv}_{k,1}$, then it follows that
$\varphi = g_2^{-1} p \, g_1^{-1}$ can also be expressed as a
product over a fixed finite set of elements of ${\it Inv}_{k,1}$.
We assume that a finite generating set for $G_{k,1}$ has been chosen.

\smallskip

For any element $\varphi \in {\it Inv}_{k,1}$ with domain code
domC$(\varphi) = P$ and image code imC$(\varphi) = Q$, we distinguish four
cases, depending on
the maximality or non-maximality of $P$ and $Q$.

\smallskip

\noindent (1) If $P$ and $Q$ are both maximal prefix codes then
$\varphi \in G_{k,1}$, and  we can express $\varphi$ over a finite fixed
generating set of $G_{k,1}$.

\smallskip

\noindent (2) Assume $P$ and $Q$ are both non-maximal prefix codes. By
Lemma \ref{ext_to_max_pref_code} there are finite maximal prefix codes
$P',Q'$ such that $P \subset P'$, $Q \subset Q'$, and $|P'| = |Q'|$; and
by Lemma \ref{card_max_pref_code}, $|P'| = |Q'| = 1 + (k-1)N$ for some
$N \geq 0$. Consider the following maximal prefix code $C$, of cardinality
$|P'| = |Q'| = 1 + (k-1)N$:

\medskip

 \ \ \  
 $C \ = \ \bigcup_{r=0}^{N-2} a_1^r (A-\{a_1\}) \ \ \cup \ a_1^{N-1} A$.

\medskip

\noindent The maximal prefix code $C$ is none other than the code $Q_{i,j}$
when $i=k$ and $j=N-1$ \, (introduced in the proof of Lemma
\ref{special_pref_code}, Fig.\ 1). 
The elements $g_1: C \to P'$ and $g_2: Q' \to C$ of $G_{k,1}$
can be chosen so that $\psi = g_2 \varphi g_1(.)$ is a partial identity
with ${\rm domC}(\psi) = {\rm imC}(\psi) \subset C$ consisting of the $|P|$ 
first elements of $C$ in the dictionary order. So, $\psi$ is the identity map
restricted to these $|P|$ first elements of $C$, and $\psi$ is undefined on
the rest of $C$. To describe ${\rm domC}(\psi) = {\rm imC}(\psi)$ in more
detail, let us write $|P| = i + (k-1) \, \ell$, for some 
$i, \ell$ with $1 \leq i < k$ and $0 \leq \ell \leq N-1$. Then

\medskip

 \ \ \
 ${\rm domC}(\psi) = {\rm imC}(\psi) \ = \ $
 $a_1^{N-1} A \ \cup \ \bigcup_{r=j+1}^{N-2} a_1^r (A-\{a_1\}) \ \ \cup \ $ 
 $a_1^j \ \{a_2, \ldots, a_i\}$. 

\medskip

\noindent where $j = N-1-\ell$. Since \ $\psi = {\rm id}_{{\rm domC}(\psi)}$, 
we claim:

\smallskip

\noindent By essential maximal extension \

\smallskip

$\psi \ = \ {\rm id}_S$ \ (as elements of ${\it Inv}_{k,1}$), \ \ where \
      $S = {\{a_1^j a_1, a_1^j a_2, \ldots, a_1^j a_i\}}$,

\smallskip

\noindent with $i, j$ as in the description of ${\rm domC}(\psi) = {\rm
imC}(\psi)$ above, i.e., $1 < i < k$, \ $N-1 \geq j = N-1-\ell \geq 0$, and 
$|P| = i + (k-1) \, \ell$.

Indeed, if $|P| < k$ then $S$ is just ${\rm domC}(\psi)$, with $i = |P|$,
and $\ell=0$ (hence $j = N-1$).  
If $|P| \geq k$ then the maximum essential extension of $\psi$ will replace 
the $1+ (k-1) \, \ell$ elements 
 \ $a_1^{N-1} A \ \cup \ \ \bigcup_{r=N-j+1}^{N-2} a_1^r (A-\{a_1\})$ \
by the single element $a_1^{N-\ell+1} = a_1^{j+1}$. 
What remains is the set \

\smallskip

$S \ = \ \{a_1^{j+1}\} \ \cup \ a_1^j \, \{a_2, \ldots, a_i\}$.

\smallskip

Finally, by Lemma \ref{special_generation}, ${\rm id}_S$ (where
$S = {\{a_1^j a_1, a_1^j a_2, \ldots, a_1^j a_i\}}$)
can be generated by the $k+1$ elements \
 \ $\{(a_1 \mapsto a_1a_1), \ (a_1a_1 \mapsto a_1)\} \ \cup \ $ 
$\{ {\rm id}_{\{a_1a_1, \ a_1a_2, \ \ldots, \ a_1a_i\}} : $
$ 1 \leq i \leq k-1 \}$.

\smallskip

\noindent (3) Assume $P$ is a maximal prefix code and $Q$ is non-maximal.
Let $Q'$ be the finite maximal prefix code obtained by saturating the prefix
tree of $Q$. Then $Q \subset Q'$, \ $|Q'| = 1 + (k-1)N'$, and
$|P| = 1 + (k-1)N$ for some $N' > N \geq 0$. We consider the maximal prefix
codes $C$ and $C'$ as defined in the proof of (2), using $N'$ for defining 
$C'$. We can choose $g_1: C \to P$ and $g_2: Q' \to C'$ in $G_{k,1}$ so that 
$\psi = g_2 \varphi g_1(.)$ is the dictionary-order preserving map that
maps $C$ to the first $|C|$ elements of $C'$. So we have  

\smallskip

${\rm domC}(\psi) = C$, \ and  

\smallskip

${\rm imC}(\psi) = S_0$ , where $S_0 \subset C'$ consist of the 
$|C|$ first elements of $C'$, in dictionary order. 

\smallskip

\noindent Since $|C| = 1 + (k-1) \, N$, we can describe $S_0$ in more detail 
by 

\smallskip

$S_0 \ = \ \bigcup_{r=N'-N}^{N'-2} a_1^r (A-\{a_1\}) \ \ \cup \ a_1^{N'-1} A$.

\smallskip

\noindent Next, by essential maximal extension we now obtain 
 \ $\psi \ = \ (\varepsilon \mapsto a_1^{N'-N})$.

Indeed, we saw that  $|P| = 1 + (k-1) \, N$. If $|P| = 1$ then 
$P = \{\varepsilon\}$, and $\psi \ = \ (\varepsilon \mapsto a_1^{N'})$.
If $|P| \geq k$ then maximum essential extension of $\psi$ will replace
all the elements of $C$ by the single element $\varepsilon$, and it will
replace all the elements of $S_0$ by the single element $a_1^{N'-N}$.

Finally, by Lemma \ref{special_generation}, 
$(\varepsilon \mapsto a_1^{N'-N})$ is generated by the two elements 
$(\varepsilon \mapsto a_1)$ and $(a_1 \mapsto a_1 a_1)$.

\smallskip

\noindent (4) The case where $P$ is a non-maximal maximal prefix code and $Q$
is maximal can be derived from case (3) by taking the inverses of the
elements from case (3).
 \ \ \ $\Box$


\begin{thm} \label{M_fin_gen} \ \
The monoid $M_{k,1}$ is finitely generated.
\end{thm}
{\bf Proof.} Let $\varphi: P \to Q$ be the table of any element of $M_{k,1}$,
mapping $P$ onto $Q$, where $P, Q \subset A^*$ are finite prefix codes.
The map described by the table is total and surjective, so if
 $|P| = |Q|$ (and in particular, if $\varphi$ is the empty map) then
$\varphi \in {\it Inv}_{k,1}$, hence $\varphi$ can be
expressed over the finite generating set of ${\it Inv}_{k,1}$.
In the rest of the proof we assume $|P| > |Q|$. The main observation is the
following.

\smallskip

\noindent
{\sf Claim.} \ $\varphi$ can be written as the composition of finitely many
elements $\varphi_i \in M_{k,1}$
with tables $P_i \to Q_i$ such that $0 \leq |P_i| - |Q_i| \leq 1$.

\smallskip

\noindent
{\sf Proof of the Claim:}
We use induction on $|P| - |Q|$. There is nothing to prove when
$|P| - |Q| \leq 1$, so we assume now that $|P| - |Q| \geq 2$.

If $\varphi(x_1) = \varphi(x_2) =  \varphi(x_3) = y_1$ for some
$x_1, x_2, x_3 \in P$ (all three being different) and $y_1 \in Q$, then we 
can write $\varphi$ as a composition $\varphi(.) = \psi_2 \circ \psi_1(.)$, 
as follows.  The map $\psi_1: P \longrightarrow P - \{x_1\}$ is defined by
$\psi_1(x_1) = \psi_1(x_2) = x_2$, and acts as the identity
everywhere else on $P$. The map $\psi_2: P - \{x_1\} \longrightarrow Q$ is 
defined by $\psi_2(x_2) = \psi_2(x_3) = y_1$, and acts in the same way as
$\varphi$ everywhere else on $P - \{x_1\}$.
Then for $\psi_1$ we have \ \ $|P| -|P - \{x_1\}| \ < \ |P| - |Q|$, \ and 
for $\psi_2$ we have \ \ $|P - \{x_1\}| - |Q| \ < \ |P| - |Q|$.

If $\varphi(x_1) = \varphi(x_2) =  y_1$ and
$\varphi(x_3) = \varphi(x_4) = y_2$ for some $x_1, x_2, x_3, x_4 \in P$
(all four being different) and $y_1, y_2 \in Q$ ($y_1 \neq y_2$), then we
can write $\varphi$ as a composition $\varphi(.) = \psi_2 \circ \psi_1(.)$, 
as follows.  
First the map $\psi_1: P \longrightarrow P - \{x_1\}$ is defined by
$\psi_1(x_1) = \psi_1(x_2) = x_2$, and acts as the identity everywhere else 
on $P$. Second, the map $\psi_2: P - \{x_1\} \longrightarrow Q$ is defined
by $\psi_2(x_2) = y_1$ and $\psi_2(x_3) = \psi_2(x_4) = y_2$, and acts like
$\varphi$ everywhere else on $P - \{x_1\}$. Again, for $\psi_1$ we have  
 \ \   $|P| -|P - \{x_1\}| \ < \ |P| - |Q|$ \ and for $\psi_2$ we have 
 \ \  $|P - \{x_1\}| - |Q| \ < \ |P| - |Q|$.
 \ \ \ \ \ [End, proof of the Claim.]

\medskip

Because of the Claim we now only need to consider elements 
$\varphi \in M_{k,1}$ with tables $P \to Q$ such that the prefix codes 
$P, Q$ satisfy $|P| = |Q| + 1$. We denote $P = \{p_1, \ldots, p_n\}$ and
$Q = \{ q_1, \ldots, q_{n-1}\}$, with $\varphi(p_j) = q_j$ for 
$1 \leq j \leq n-1$, and $\varphi(p_{n-1}) = \varphi(p_n) = q_{n-1}$.
We define the following prefix code $C$ with $|C| = |P|$:

\smallskip

\noindent $\bullet$ \ if $|P| = i \leq k$ 
then \ \ $C \ = \ \{a_1, \ldots, a_i \}$; \ note 
that $i \geq 2$, since $|P| > |Q| >0$;

\smallskip

\noindent $\bullet$ \ if $|P| > k$ then \ \
$C \ = \ \{a_2, \ldots, a_i \} \ \cup \ $
       $\bigcup_{r=1}^{j-1} a_1^r (A-\{a_1\}) \ \ \cup \ \ a_1^j A$,

\smallskip

\noindent where $i,j$ are such that $|P| = i + (k-1)j$, \ $2 \leq i \leq k$,
and $1 \leq j$ (see Fig.\ 1).
Let us write $C$ in increasing dictionary order as
$C = \{c_1, \ldots, c_n\}$. The last element of $C$ in the dictionary
order is thus $c_n = a_i$.

We now write $\varphi(.) = \psi_3 \, \psi_2 \, \psi_1(.)$ where
$\psi_1$, $\psi_2$, $\psi_3$ are as follows: \\
$\bullet$ \ $\psi_1: P \longrightarrow C$ \ is bijective and is defined by
 $p_j \mapsto c_j$ for $1 \leq j \leq n$; \\
$\bullet$ \ $\psi_2: C \longrightarrow C - \{ a_i\}$ \ is the identity map on
$\{c_1, \ldots, c_{n-1} \}$, and $\psi_2(c_n) = c_{n-1}$. \\
$\bullet$ \ $\psi_3: C - \{ a_i\} \longrightarrow Q$ \ is bijective and is 
defined by $c_j \mapsto q_j$ for $1 \leq j \leq n-1$.

It follows that $\psi_1$ and $\psi_3$ can be expressed over the finite
generating set of ${\it Inv}_{k,1}$.
On the other hand, $\psi_2$ has a maximum essential extension, as follows.

\medskip

\noindent
$\bullet$ \ If \ $2 \leq |P| = i \leq k$ \ then

\medskip

$\psi_2 \ = \ $
$\left( \hspace{-.08in} \begin{array}{l|l|l|l|l}
a_1 & \ \ldots \ & a_{i-2} & a_{i-1} & a_i  \\
a_1 & \ \ldots \ & a_{i-2} & a_{i-1} & a_{i-1}
\end{array} \hspace{-.08in} \right) \ = \ $
$\left( \hspace{-.08in} \begin{array}{l|l}
{\rm id}_{\{a_1, \ \ldots, \ a_{i-1}\}}  & a_i  \\
                                  \      & a_{i-1}
\end{array} \hspace{-.08in} \right)$.

\medskip

\noindent
$\bullet$ \ If \ $|P| = i + (k-1)j > k$ \ and if $i > 2$ then, after
maximal essential extension, $\psi_2$ also becomes \ 

\medskip

$\max(\psi_2) \ = \ $
$\left( \hspace{-.08in} \begin{array}{l|l}
{\rm id}_{\{a_1, \ \ldots, \ a_{i-1}\}}  & a_i  \\
                                  \      & a_{i-1}
\end{array} \hspace{-.08in} \right)$. 

\medskip

\noindent 
$\bullet$ \ If \ $|P| = i + (k-1)j > k$ \ and if $i = 2$ then, after
essential extensions,

\medskip

$\max(\psi_2) \ = \ $
$\left( \hspace{-.08in} \begin{array}{l|l|l|l|l|l}
a_1a_1 & \ \ldots \ & a_1 a_{k-2} & a_1a_{k-1} & a_1a_k & a_2 \\
a_1a_1 & \ \ldots \ & a_1 a_{k-2} & a_1a_{k-1} & a_1a_k & a_1a_k
\end{array} \hspace{-.08in} \right) \ = \ $
$\left( \hspace{-.08in} \begin{array}{l|l}
{\rm id}_{a_1A}  & a_2  \\
          \      & a_1a_k
\end{array} \hspace{-.08in} \right) \ = \ $
$\left( \hspace{-.08in} \begin{array}{l|l}
a_1 & a_2  \\
a_1 & a_1a_k
\end{array} \hspace{-.08in} \right) $.

\medskip

\noindent
In summary, we have factored $\varphi$ over a finite set of generators of
${\it Inv}_{k,1}$ and $k$ additional generators in $M_{k,1}$.
 \ \ \ $\Box$

\bigskip

\noindent {\bf Factorization algorithm:} The proofs of Theorems 
\ref{Inv_fin_gen} and \ref{M_fin_gen} are constructive; they provide 
algorithms that, given $\varphi \in {\it Inv}_{k,1}$ or $ \in M_{k,1}$,
output a factorization of  $\varphi$ over the finite generating set of
${\it Inv}_{k,1}$, respectively $M_{k,1}$.

\bigskip

In \cite{Hig74} (p.\ 49) Higman introduces a four-element generating set for
$G_{2,1}$; a special property of these generators is that their domain codes
and their image codes only contain words of length $\leq 2$, and that
$\big| \, |\gamma(x)| - |x| \, \big| \leq 1$ for every generator $\gamma$ 
and every $x \in {\rm domC}(\gamma)$. The generators in the finite generating 
set of $M_{k,1}$ that we introduced above also have those properties. Thus we 
obtain:

\begin{cor} \label{lengths_in_domCimC} \
The monoid $M_{2,1}$ has a finite generating set such that all the generators
have the following property: The domain codes and the image codes only contain
words of length $\leq 2$, and \ $\big| \, |\gamma(x)| - |x| \, \big| \leq 1$
for every generator $\gamma$ and every $x \in {\rm domC}(\gamma)$.
 \ \ \ \ \ $\Box$
\end{cor}

\noindent
For reference we list an explicit {\em finite generating set for} $M_{2,1}$. 
It consists, first, of the Higman generators of $G_{2,1}$ (\cite{Hig74} 
p.\ 49):

\medskip

{\sc Not} $=$ $\left( \hspace{-.08in} \begin{array}{l|l}
0 & 1  \\
1 & 0
\end{array} \hspace{-.08in} \right)$, \ \
$(01 \leftrightarrow 1) = $
$\left( \hspace{-.08in} \begin{array}{r|r|r}
00 & 01 & 1  \\
00 &  1 & 01
\end{array} \hspace{-.08in} \right)$, \ \
$(0 \leftrightarrow 10) = $
$\left( \hspace{-.08in} \begin{array}{r|r|r}
0 & 10 & 11  \\
10 & 0 & 11
\end{array} \hspace{-.08in} \right)$, \ \ and 

\smallskip

$\tau_{1,2} = $
$\left( \hspace{-.08in} \begin{array}{r|r|r|r}
00 & 01 & 10 & 11  \\
00 & 10 & 01 & 11
\end{array} \hspace{-.08in} \right)$;

\medskip

\noindent the additional generators for ${\it Inv}_{2,1}$:

\medskip

$(\varepsilon \to 0)$, \ \ $(0 \to \varepsilon)$, \ \ $(0 \to 00)$, \ \
$(00 \to 0)$;

\medskip

\noindent the additional generators for $M_{2,1}$:

\medskip

$\left( \hspace{-.08in} \begin{array}{l|l}
0 & 1  \\
0 & 0
\end{array} \hspace{-.08in} \right)$, \ and \ \
$\left( \hspace{-.08in} \begin{array}{r|r}
0 & 1  \\
0 & 01
\end{array} \hspace{-.08in} \right)$
$ = \left( \hspace{-.08in} \begin{array}{r|r|r}
00 & 01 & 1 \\
00 & 01 & 01
\end{array} \hspace{-.08in} \right)$ .

\bigskip

\noindent Observe that Higman's generators of $G_{k,1}$ (in \cite{Hig74}
p.\ 27) have domain and image codes with at most 3 internal
vertices. We observe that the additional generators that we introduced for
${\it Inv}_{k,1}$ and $M_{k,1}$ have domain and image codes have at most 2
internal vertices.

\bigskip

\noindent The following problem remains open: Are ${\it Inv}_{k,1}$ and 
$M_{k,1}$ finitely {\bf presented}?


\section{The word problem of the Thompson-Higman monoids}

We saw that the Thompson-Higman monoid $M_{k,1}$ is finitely generated.
We want to show now that the word problem of $M_{k,1}$ over any 
finite generating set can be decided in deterministic polynomial time, i.e.,
it belongs to the complexity class {\sf P}. 
\footnote{\, This section has been revised in depth, to correct errors.}

In \cite{BiThomps} it was shown that the word problem of the Thompson-Higman
group $G_{k,1}$ over any finite generating set is in {\sf P}. In fact, it
is in the parallel complexity class {\sf AC}$_1$ \cite{BiThomps}, and it is
co-context-free \cite{LehSchw}.
In \cite{BiCoNP} it was shown that the word problem of the Thompson-Higman
group $G_{k,1}$ over the infinite generating set \
$\Gamma_{k,1} \cup \{ \tau_{i,i+1} : i > 0\}$ \ is {\sf coNP}-complete,
where $\Gamma_{k,1}$ is any finite generating set of $G_{k,1}$;
the position transposition $\tau_{i,i+1} \in G_{k,1}$ has \
${\rm domC}(\tau_{i,i+1}) = {\rm imC}(\tau_{i,i+1}) = A^{i+1}$, and is
defined by \ $u \alpha \beta \mapsto u \beta \alpha$ \ for all letters
$\alpha, \beta \in A$ and all words $u \in A^{i-1}$.
We will see below that the word problem of $M_{k,1}$ over
$\Gamma_{k,1} \cup \{ \tau_{i,i+1} : i > 0\}$ is also {\sf coNP}-complete,
where $\Gamma_{k,1}$ is any finite generating set of $M_{k,1}$.


\subsection{The image code formula} 

Our proof in \cite{BiThomps} that the word problem of the Thompson-Higman
group $G_{k,1}$ (over any finite generating set) is in {\sf P}, was based on 
the following fact (the {\em table size formula}):

\medskip

 \ \ \ \ \  $\forall \varphi, \psi \in G_{k,1}$: 
 \ \ $\|\psi \circ \varphi\| \leq \|\psi\| + \|\varphi\|$.

\medskip 

\noindent Here $\|\varphi\|$ denotes the {\em table size} of $\varphi$, 
i.e., the cardinality of ${\rm domC}(\varphi)$. See Proposition 3.5, 
Theorem 3.8, and Proposition 4.2 in \cite{BiThomps}. In $M_{k,1}$ the above 
formula does not hold in general, as the following example shows. 
We give some definitions and notation first.

\begin{defn}
For any finite set $S \subseteq A^*$ we denote the length of the longest
word in $S$ by $\ell(S)$.
The cardinality of $S$ is denoted by $|S|$.

The {\em table} of a right-ideal morphism $\varphi$ is the set
$\{(x, \varphi(x)) : x \in {\rm domC}(\varphi)\}$.
\end{defn}

\begin{pro} \label{exp_table} \    For every $n>0$ there exists
$\Phi_n = \varphi_2^{n-1} \varphi_1 \in M_{2,1}$ (for some
$\varphi_1, \varphi_2 \in  M_{2,1}$) with the following properties:

The table sizes are $\|\Phi_n\| = 2^n$, and
 \ $\|\varphi_2\| = \|\varphi_1\| = 2$.
So, $\|\Phi_n\|$ is exponentially larger than
 \, $(n-1) \cdot \|\varphi_2\| + \|\varphi_1\|$.
Hence the table size formula does not hold in $M_{2,1}$.

The word lengths of $\varphi_1, \varphi_2$, and
$\Phi_n$ (over the finite generating set $\Gamma$ of $M_{2,1}$ from 
Section 3 in \cite{BiMonTh}) satisfy
 \ $|\varphi_1|_{_{\Gamma}} =1$, \ $|\varphi_2|_{_{\Gamma}} \leq 2$, and 
$|\Phi_n|_{_{\Gamma}} < 2 n$.
So the table size of $\Phi_n$ is exponentially larger than its word length:
 \ $\|\Phi_n\| \, > \, \sqrt{2}^{\, |\Phi_n|_{_{\Gamma}}}$.
\end{pro}
{\bf Proof.} Consider $\varphi_1, \varphi_2 \in M_{2,1}$ given by the tables
 \ $\varphi_1 = \{(0 \mapsto 0), \ (1 \mapsto 0)\}$, and
 \ $\varphi_2 = \{(0 0 \mapsto 0), \ (0 1 \mapsto 0)\}$.
One verifies that \ $\Phi_n = \varphi_2^{n-1} \circ \varphi_1(.)$ \ sends
every bitstring of length $n$ to the word $0$; its domain code
is $\{0,1\}^n$, its image code is $\{0\}$, and it is its maximum essential
extension. Thus, $\|\varphi_2^{n-1} \circ \varphi_1\| \ = \ 2^n$,
whereas \ $(n-1) \cdot \|\varphi_2\| + \|\varphi_1\| = 2 n$.
Also, \ $\varphi_2(.) \, = \, $
$(0 \mapsto 0, 1 \mapsto 0) \cdot (0 \mapsto \varepsilon)$, so
$|\varphi_1|_{_{\Gamma}} = 1$, \ $|\varphi_2|_{_{\Gamma}} \leq 2$, and 
$|\Phi_n|_{_{\Gamma}} \leq 2n-1$; hence 
$\|\Phi_n\| > 2^{|\Phi_n|_{_{\Gamma}}/2}$. 
 \ \ \ $\Box$

\bigskip

We will use the following facts that are easy to prove:
If $R \subset A^*$ is a right ideal and $\varphi$ is a right-ideal morphism 
then $\varphi(R)$ and $\varphi^{-1}(R)$ are right ideals. 
The intersection and the union of right ideals are right ideals.
We also need the following result.

\medskip

\noindent {\bf (Lemma 3.3 of \cite{BiThomps})}
 \ {\em If $P,Q,S \subseteq A^*$ are such that $PA^* \cap QA^* = SA^*$, 
and if $S$ is a prefix code then $S \subseteq P \cup Q$.  }
 \ \ \ $\Box$

\begin{lem} \label{Im_of_rightideal} \  
Let $\theta$ be a right-ideal morphism, and assume 
$SA^* \subseteq {\rm Dom}(\theta)$, where $S \subset A^*$ is a finite prefix 
code.
Then there is a finite prefix code $R \subset A^*$ such that 
$\, \theta(SA^*) = RA^* \,$ and $\, R \subseteq \theta(S)$.
\end{lem}
{\bf Proof.} Since $\theta$ is a right-ideal morphism we have \  
$\theta(SA^*) = \theta(S) \ A^*$. Since $\theta(S)$ might not be a prefix
code we take \ $R = \{r \in \theta(S) : r$ is minimal (shortest) in the
prefix order within $\theta(S)\}$. Then $R$ is a prefix code that has the
required properties.
  \ \ \ $\Box$

\begin{lem} \label{InvIm_of_rightideal}\footnote{\, This Lemma was incorrect 
in the earlier versions of this paper and in \cite{BiMonTh}.}
 \ For any right-ideal morphism $\theta$ and any prefix code $Z \subset A^*$,
$\theta^{-1}(Z)$ is a prefix code. 

In particular, $\theta^{-1}({\rm imC}(\theta))$ is a prefix code, and 
$\, \theta^{-1}({\rm imC}(\theta)) \subseteq {\rm domC}(\theta)$.
There exist right-ideal morphisms $\theta$ with finite domain code, such 
that $\, \theta^{-1}({\rm imC}(\theta)) \neq {\rm domC}(\theta)$. 
\end{lem}
{\bf Proof.}  First, $\theta^{-1}(Z)$ is a prefix code. Indeed, if we had
$x_1 = x_2u$ for some $x_1, x_2 \in \theta^{-1}(Z)$ with $u$ non-empty, then
$\theta(x_1) = \theta(x_2) \ u$, with $\theta(x_1), \theta(x_2) \in Z$. This 
would contradict the assumption that $Z$ is a prefix code.

Second, let $Q = {\rm imC}(\theta)$; then
$\theta^{-1}(Q) \, A^* \subseteq \theta^{-1}(QA^*)$. 
Indeed, if $x \in \theta^{-1}(Q)$, then $x = pw$ for some 
$p \in {\rm domC}(\theta)$ and $w \in A^*$. Hence,
$\theta(x) = \theta(p) \, w$, and $\theta(x) \in Q$. 
Since $\theta(p) \, w \in Q$ and $\theta(p) \in Q A^*$, we have 
$\theta(p) \, w = \theta(p)$ (since $Q$ is a prefix code). So $w$ is 
empty, hence $x = pw = p \in {\rm domC}(\theta)$.  

Example: Let $A = \{0,1\}$, and let $\theta$ be the right-ideal morphism 
defined by ${\rm domC}(\theta) = \{01, 1\}$, 
${\rm imC}(\theta) = \{\varepsilon\}$, and
 $\, \theta(01) = 0$, $\, \theta(1) = \varepsilon$.
Then, $\theta^{-1}({\rm imC}(\theta)) = \{1\} \neq {\rm domC}(\theta)$.
  \ \ \ $\Box$

\medskip

\noindent The following generalizes the ``table size formula'' of 
$G_{k,1}$ to the monoid $M_{k,1}$.

\begin{thm} \label{sum_of_imC} {\bf (Generalized image code formulas).}
\footnote{\, This Theorem was incorrect in the previous versions and in
\cite{BiMonTh}; this is a corrected (and expanded) version.} \\   
Let $\varphi_i$ be right-ideal morphism with finite domain codes, for 
$i = 1, 2, \ldots, n$.  Then

\bigskip

\noindent {\bf (1)} \ \ \ \   
$\big|{\rm imC}(\varphi_n \circ \ \ldots \ \circ \varphi_1)\big|$
$ \ \leq \ $    
$|{\rm imC}(\varphi_1)|  \ + \ $
$\sum_{i=2}^n |\varphi_i({\rm domC}(\varphi_i))|$,  

\bigskip

\noindent {\bf (2)} \ \ \ \    
$\ell \big({\rm domC}(\varphi_n \circ \ \ldots \ \circ \varphi_1)\big)$ 
$ \ \leq \ $
$\sum_{i=1}^n \ell ({\rm domC}(\varphi_i))$,

\bigskip

\noindent {\bf (3)} \ \ \ \    
$\ell \big(\varphi_n \ldots \varphi_1({\rm domC}(\varphi_n
 \circ \ \ldots \ \circ \varphi_1)) \big)$
$ \ \leq \ $
$\sum_{i=1}^n \ell(\varphi_i({\rm domC}(\varphi_i)))$, 

\bigskip

\noindent {\bf (4)} \ \ \ \
$\ell \big({\rm imC}(\varphi_n \circ \ \ldots \ \circ \varphi_1) \big)$
$ \ \leq \ $    
$\ell({\rm imC}(\varphi_1)) \ + \ $
$\sum_{i=2}^n \ell(\varphi_i({\rm domC}(\varphi_i)))$,

\bigskip

\noindent {\bf (5)} \ \ \ \
$\big|\varphi_n \ldots \varphi_1({\rm domC}(\varphi_n
 \circ \ \ldots \ \circ \varphi_1)) \big|$
$ \ \leq \ $
$\sum_{i=1}^n |(\varphi_i({\rm domC}(\varphi_i))|$, \ and

\bigskip

\hspace{.25in} 
$\varphi_n \ldots \varphi_1({\rm domC}(\varphi_n \circ \ \ldots \, $
$\circ$ $\varphi_1))$ 
$ \ \subseteq \ \ $ $\bigcup_{i=1}^n$
$\varphi_n \ldots \varphi_i({\rm domC}(\varphi_i))$.
\end{thm}
{\bf Proof.} \ Let $P_i = {\rm domC}(\varphi_i)$ and 
$Q_i = {\rm imC}(\varphi_i)$.

\smallskip

\noindent (1) \  The proof is similar to the proof of Proposition 3.5
in \cite{BiThomps}. We have \  ${\rm Dom}(\varphi_2 \circ \varphi_1) = $
$\varphi_1^{-1}(Q_1A^* \cap P_2A^*)$ \ and 
 \ ${\rm Im}(\varphi_2 \circ \varphi_1) = \varphi_2(Q_1A^* \cap P_2A^*)$. 
So the following maps are total and onto on the indicated sets:

 \ \ \  \ \ \    
$\varphi_1^{-1}(Q_1A^* \cap P_2A^*) \ \stackrel{\varphi_1}{\longrightarrow}$
$ \ Q_1A^* \cap P_2A^* \ \stackrel{\varphi_2}{\longrightarrow} \ $
$\varphi_2(Q_1A^* \cap P_2A^*)$.

\smallskip

\noindent By Lemma 3.3 of \cite{BiThomps} (quoted above) we have
$Q_1A^* \cap P_2A^* = SA^*$ for some finite prefix code $S$ with 
$S \subseteq Q_1 \cup P_2$. Moreover, by Lemma \ref{Im_of_rightideal} we 
have $\varphi_2(SA^*) = R_2A^*$ for some finite prefix code $R_2$ such that 
$R_2 \subseteq \varphi_2(S)$. Now, since $S \subseteq Q_1 \cup P_2$ we have 
 \ $R_2 \subseteq \varphi_2(S) \subseteq $
$\varphi_2(Q_1) \, \cup \, \varphi_2(P_2)$. 
Thus, \ $|{\rm imC}(\varphi_2 \circ \varphi_1)| = $
$|R_2| \leq |\varphi_2(P_2)| + |\varphi_2(Q_1)|$. Since 
$|\varphi_2(Q_1)| \le |Q_1|$, we have 
$|R_2| \leq |\varphi_2(P_2)| + |Q_1|$.

By induction for $n > 2$, $|{\rm imC}(\varphi_n$
$\circ$ $\varphi_{n-1} \circ \ \ldots \ \circ \varphi_1)|$  $\le$
$|\varphi_n({\rm domC}(\varphi_n))|$ $+$
$|{\rm imC}(\varphi_{n-1} \circ \ \ldots \ \circ \varphi_1)|$ $\le$
$|\varphi_n({\rm domC}(\varphi_n))|$ $+$
$\sum_{i=2}^{n-1} |\varphi_i({\rm domC}(\varphi_i))| \ + \
|{\rm imC}(\varphi_1)|$ .

\medskip

\noindent (2) \ We prove the formula when $n = 2$; the general formula
then follows immediately by induction.
Let $x \in {\rm domC}(\varphi_2 \circ \varphi_1)$; then $\varphi_1(x)$ is
defined, hence $x = p_1 u$ for some $p_1 \in P_1$, $u \in A^*$.  And 
$\varphi_2$ is defined on $\varphi_1(x) = \varphi_1(p_1) \, u$, so 
$\varphi_1(x) \in P_2 A^* = {\rm Dom}(\varphi_2)$.  Hence there exist
$p_2 \in P_2$ and $v \in A^*$ such that 

\medskip

\noindent $(\star)$ \hspace{.5in}  $\varphi_1(p_1) \, u = p_2 v$ 
 \ $\in \varphi_1(P_1) \, A^* \cap P_2 A^*$.  

\medskip

\noindent It follows that $u$ and $v$ are suffix-comparable.  

\medskip

\noindent {\sf Claim.} \ {\it The words $u$ and $v$ in $(\star)$ satisfy:
$\, u = \varepsilon$, or $v = \varepsilon$.    } 

\smallskip

\noindent Proof of the Claim: Since $u$ and $v$ are suffix-comparable, let 
us first consider the case where $v$ is a suffix of $u$, i.e., $u = t v$ for 
some $t \in A^*$. Then $\varphi_1(x) = \varphi_1(p_1) \, t v = p_2 v$, hence 
$\varphi_1(p_1) \, t = p_2$, hence 
$\varphi_2$ is defined on $\varphi_1(p_1) \, t = p_2$. So, 
$\varphi_2 \circ \varphi_1$ is defined on $p_1 t$, i.e., 
$p_1 t \in {\rm domC}(\varphi_2 \circ \varphi_1)$. 
But we also have $x = p_1 t v \in {\rm domC}(\varphi_2 \circ \varphi_1)$.
Since ${\rm domC}(\varphi_2 \circ \varphi_1)$ is a prefix code, it follows 
that $v = \varepsilon$.

Let us next consider the other case, namely where $u$ is a suffix of $v$, 
i.e., $v = s u$ for some $s \in A^*$. Then 
$\varphi_1(x) = \varphi_1(p_1) \, u = p_2 s u$, hence 
$\varphi_1(p_1) = p_2 s$, hence 
$\varphi_2$ is defined on $\varphi_1(p_1) = p_2 s$, hence 
$p_1 \in {\rm domC}(\varphi_2 \circ \varphi_1)$. But we also have 
$x = p_1 u \in {\rm domC}(\varphi_2 \circ \varphi_1)$.
Since ${\rm domC}(\varphi_2 \circ \varphi_1)$ is a prefix code, it follows
that $u = \varepsilon$.  \ \ \ [This proves the Claim.]

\medskip

Now for $x \in {\rm domC}(\varphi_2 \circ \varphi_1)$ we have
$x = p_1 u$, and $\varphi_1(p_1) \, u = p_2 v$, hence $|x| = |p_1| + |u|$
and $|\varphi_1(p_1)| + |u| = |p_2| + |v|$. 
By the Claim, either $|u| = 0$ or $|v| = 0$.

If $\, |u| = 0 \,$ then $\, |x| = |p_1| \le \ell({\rm domC}(\varphi_1))$. 

If $\, |v| = 0 \,$ then $\, |x| = |p_1| + |u|$ $=$ 
$|p_1| + |p_2| + |v| - |\varphi_1(p_1)|$ $=$ 
$|p_1| + |p_2| - |\varphi_1(p_1)| \le |p_1| + |p_2|$ $\le$
$\ell({\rm domC}(\varphi_1)) + \ell({\rm domC}(\varphi_2))$.

\medskip

\noindent (3) \ As in the proof of (2) we only need to consider $n = 2$.  
Let $x \in$ ${\rm domC}(\varphi_2 \varphi_1)$, hence 
$\varphi_2 \varphi_1(x) \in$ 
$\varphi_2 \varphi_1({\rm domC}(\varphi_2 \varphi_1))$.  By $(\star)$ (and 
with the notation of the proof of (2)) we have 
$\varphi_2 \varphi_1(x) = \varphi_2(\varphi_1(p_1) \, u)$ $=$ 
$\varphi_2(p_2) \, v \in \varphi_2(\varphi_1(P_1) \, A^* \cap P_2 A^*) =$
${\rm Im}(\varphi_2 \varphi_1)$.  By the reasoning of the proof of (2), 
we have two cases:

If $|u| = 0$ then 
$|v| = |\varphi_1(p_1)| + |u| - |p_2| = |\varphi_1(p_1)| - |p_2|$ $\le$
$|\varphi_1(p_1)|$.
Hence, $|\varphi_2 \varphi_1(x)| = |\varphi_2(p_2)| + |v| \le$
$|\varphi_2(p_2)| + |\varphi_1(p_1)| \ \le \ $
$\ell(\varphi_2({\rm domC}(\varphi_2))$ $+$ 
$\ell(\varphi_1({\rm domC}(\varphi_1))$.

If $|v| = 0$ then $\varphi_2 \varphi_1(x) = \varphi_2(p_2)$, hence
 \ $|\varphi_2 \varphi_1(x)| = |\varphi_2(p_2)| \le $
$\ell(\varphi_2({\rm domC}(\varphi_2))$.

\medskip

\noindent (4) \ We first consider the case $n=2$. As we saw in the proof 
of (1), ${\rm imC}(\varphi_2 \varphi_1) = R_2$ where 
$R_2 \subseteq \varphi_2(S)$, and where $S$ is a prefix code such that 
$S  \subseteq Q_1 \cap P_2$. Hence $R_2 \subseteq$ 
$\varphi_2(Q_1) \cup \varphi_2(P_2)$. 
 
Hence for any $z \in R_2$, either $z \in \varphi_2(P_2)$ or 
$z \in \varphi_2(Q_1)$.
If $z \in \varphi_2(P_2)$ then $|z| \leq \ell(\varphi_2(P_2))$.  If
$z \in \varphi_2(Q_1)$, then $z = \varphi_2(q_1)$ for some
$q_1 \in Q_1 \cap P_2 A^*$, so $q_1 = p_2 u$ for some $p_2 \in P_2$ and 
$u \in A^*$. We have $q_1 \in P_2 A^*$ ($= {\rm Im}(\varphi_2)$), so 
$q_1 \in  {\rm Im}(\varphi_2)$.
Now $|z| = |\varphi_2(p_2)| + |u|$, and 
$|u| = |q_1| - |p_2| \le |q_1| \le \ell({\rm imC}(\varphi_1))$. 
Thus, $|z| \le |\varphi_2(p_2)| + \ell({\rm imC}(\varphi_1))$
$\le$ $\ell(\varphi_2({\rm domC}(\varphi_2)))$ $+$
$\ell({\rm imC}(\varphi_1))$.

The formula for $n > 2$ now follows by induction in the same way as in the
proof of (1).

\medskip

\noindent (5) \ We first prove the formula for $n = 2$.
As we saw in the proof of (2), if $x \in {\rm domC}(\varphi_2 \varphi_1)$ 
then there exist $u, v \in A^*$, $p_1 \in P_1$, $p_2 \in P_2$, such that 
$x = p_1 u$ and $\varphi_1(x) = \varphi_1(p_1) \, u = p_2 v$. 
Moreover, by the Claim in (2) we have $u = \varepsilon$ or 
$v = \varepsilon$.
Also, $\varphi_2 \varphi_1(x) = \varphi_2(\varphi_1(p_1) \, u) = $
$\varphi_2(p_2) \, v$.

If $v = \varepsilon$ then $\varphi_2 \varphi_1(x) = \varphi_2(p_2)$ 
$\in \varphi_2({\rm domC}(\varphi_2))$.
If $u = \varepsilon$ then $\varphi_2 \varphi_1(x) = $
$\varphi_2 \varphi_1(p_1) \in \varphi_2 \varphi_1({\rm domC}(\varphi_1))$.
Thus we proved the following fact:

\medskip

 \ \ \ \ \  $\varphi_2 \varphi_1({\rm domC}(\varphi_2 \varphi_1))$
$\, \subseteq \ $
$\varphi_2({\rm domC}(\varphi_2))$ $ \ \cup \ $
$\varphi_2 \varphi_1({\rm domC}(\varphi_1))$.

\medskip

\noindent Now, since $|\varphi_2 \varphi_1({\rm domC}(\varphi_1))| \le $
$|\varphi_1({\rm domC}(\varphi_1))|$, the fact implies that 
$\, |\varphi_2 \varphi_1({\rm domC}(\varphi_2 \varphi_1))|$ $ \ \le \ $
$|\varphi_2({\rm domC}(\varphi_2))|$ $+$
$|\varphi_1({\rm domC}(\varphi_1))|$.
By induction we immediately obtain

\medskip

$\big|\varphi_n \ldots \varphi_1({\rm domC}(\varphi_n
 \circ \ \ldots \ \circ \varphi_1)) \big|$
$ \ \leq \ $
$\sum_{i=1}^n |(\varphi_i({\rm domC}(\varphi_i))|$, \ and

\medskip

$\varphi_n \ldots \varphi_1({\rm domC}(\varphi_n \circ \ \ldots \, $
$\circ$ $\varphi_1))$

\smallskip

$ \ \subseteq \ \ $ 
$\varphi_n({\rm domC}(\varphi_n))$ $ \ \cup \ $
$\varphi_n \varphi_{n-1}({\rm domC}(\varphi_{n-1}))$
$ \ \cup \ \ \ldots \ \ \ldots \ \ \ldots$

 \ \ \  \ \ \   
$\cup \ $  $\varphi_n \ldots \varphi_i({\rm domC}(\varphi_i))$
$ \ \cup \ $  $ \ \ldots \ \ \ldots \ \ \ldots \ $ $ \ \cup \ $
$\varphi_n \ldots \varphi_i \ldots \varphi_1({\rm domC}(\varphi_1))$.
  \ \ \ \ \ $\Box$

\bigskip

\noindent {\bf Remarks.}  
Obviously, ${\rm Dom}(\varphi_2 \varphi_1)$ $\subseteq$ 
${\rm Dom}(\varphi_1)$; however, in infinitely many cases (in ``most'' 
cases), ${\rm domC}(\varphi_2 \varphi_1)$ $\not\subseteq$ 
${\rm domC}(\varphi_1)$. Instead, we have the more complicated formula of
Theorem \ref{sum_of_imC}(5).

By Prop.\ \ref{exp_table}, we cannot have a formula for 
$|{\rm domC}(\varphi_n \ldots \varphi_1)|$ of a similar nature as the 
formulas in Theorem \ref{sum_of_imC}. 

\bigskip

\noindent The following class of right-ideal morphisms plays an important
role here (as well as in Section 5 of \cite{BiCongr}, where it was 
introduced).
\footnote{\, Def.\ 4.5A, Theorem 4.5B, and Cor.\ 4.5C are new in this 
version.} 

\bigskip

\noindent {\bf Definition 4.5A (Normal).}
{\it A right-ideal morphism $\varphi$ is called {\rm normal} iff
$\, \varphi({\rm domC}(\varphi)) = {\rm imC}(\varphi)$.
} 

\bigskip

\noindent By Lemma 5.7 of \cite{BiCongr} we also have: $\, \varphi$ is 
normal iff
$\varphi^{-1}({\rm imC}(\varphi)) = {\rm domC}(\varphi)$.
In other words, $\varphi$ is normal iff $\varphi$ is entirely determined by 
the way it maps ${\rm domC}(\varphi)$ onto ${\rm imC}(\varphi)$.
 
For example, {\it every injective right-ideal morphism is normal} (by 
Lemma 5.1 in \cite{BiCongr}). The finite generating set $\Gamma$ of 
$M_{k,1}$, constructed in Section 3, consist entirely of normal right-ideal 
morphisms. 

On the other hand, the composition of two normal right-ideal morphisms does 
not always result in a normal morphism, as is shown by the following 
example: $\, {\rm domC}(f) = \{0,1\} \, $ and $\, f(0) = 0$, $f(1) = 10$;
$\, {\rm domC}(g) = \{0,1\} \, $ and $\, g(0) = g(1) = 0$; so $f$ and $g$ 
are normal. But ${\rm domC}(gf) = \{0, 1\}$ and $gf(0) = 0$, 
$\, gf(1) = 00$; so $gf$ is not normal (for more details, see Prop.\ 5.8 in 
\cite{BiCongr}).

\medskip

The next result (Theorem 4.5B) shows that every element of $M_{k,1}$ can be 
represented by a normal right-ideal morphism. So one can say informally 
that ``from the point of view of $M_{k,1}$, all right-ideal morphisms are
normal''.  For proving this we need some definitions. We always assume 
$|A| \ge 2$.  

\bigskip

\noindent {\bf Definitions and notation.} 
 \ {\it If $x_1, x_2 \in A^*$ are such that $x_1$ is a prefix of $x_2$,
i.e., $x_2 \in x_1 A^*$, we denote this by $x_1 \le_{\rm pref} x_2$.  

For $Z \subseteq A^*$, the set of {\em prefixes of $Z$} is
$\, {\sf pref}(Z) = \{v \in A^* : v \le_{\rm pref} z$ for some $z \in Z\}$.

For a set $X \subseteq A^*$ and a word $v \in A^*$, $v^{-1} X$ denotes the
set $\{s \in A^* : v s \in X\}$.

The {\rm tree of} $A^*$ has root $\varepsilon$, vertex set $A^*$,
and edge set $\, \{(w, wa) : w \in A^*, \, a \in A\}$.

A {\rm subtree} of the tree of $A^*$ has as root any string $r \in A^*$, and
as vertex set any subset $V \subseteq r A^*$, such that the following holds 
for all $v \in V$ and $u \in A^*$:
 \ $r \le_{\rm pref} u  \le_{\rm pref} v$ implies $u \in V$.
}

\medskip 

\noindent The following is a slight generalization of the classical
notion of a prefix tree.

\medskip

\noindent {\bf Definition (Prefix tree).} 
 \ {\it Let $Z \subseteq A^*$, and let $q \in A^*$. The {\rm prefix tree} 
$T(q, Z)$ is the subtree of the tree of $A^*$ with root $q$ and vertex set 
 \  $V_{q, Z} = \{v \in A^*: q \le_{\rm pref} v$, and $v \le_{\rm pref} z$ 
for some $z \in Z\}$.
}

\medskip

\noindent {\bf Remark.}
Let $L$ be the set of leaves of $T(q, Z)$; then $L$ and $q^{-1} L$ are
prefix codes.

\medskip

\noindent {\bf Definition (Saturated tree).}
 \ {\it A subtree $T$ of the tree of $A^*$ is {\rm saturated} iff for every 
vertex $v$ of $T$ we have: $v$ has no child in $T$ (i.e., $v$ is a leaf), 
or $v$ has $|A|$ children in $T$.
}

\medskip

\noindent {\bf Definition (Tree saturation).} 
\ {\it Let $T$ be a subtree of the tree of $A^*$, with root $q$, set of 
vertices $V$, and set of leaves $L$.  The {\rm saturation} of $T$ is the 
smallest (under inclusion) {\rm saturated} subtree of the tree of $A^*$ 
with root $q$, that contains $T$.  
In other words, if $T$ is just $\{q\}$, it is its own saturation; otherwise 
the saturation has root $q$ and has vertex set $\, V \cup (V - L) \cdot A$. 
We denote the saturation of $T$ by ${\rm s}T$.
} 

\medskip

\noindent {\bf Remark.} 
{\bf (1)} The prefix tree $T(q, Z)$ and its saturation have the same 
{\it depth} (i.e., length of a longest path from the root). 
Every leaf of $T(q, Z)$ is also a leaf of ${\rm s}T(q, Z)$, but  unless 
$T(q, Z)$ is already saturated, ${\rm s}T(q, Z)$ has more leaves than 
$T(q, Z)$.  
The non-leaf vertices of $T(q, Z)$ and ${\rm s}T(q, Z)$ are the same.

\noindent {\bf (2)} The number of leaves in the saturated tree ${\rm s}T(q, Z)$ is 
$\, < |V_{q, Z}| \cdot |A|$. 

\noindent {\bf (3)} Let  $L$ be the leaf set of the saturated tree
${\rm s}T(q, Z)$; 
if $Z$ is finite then $q^{-1} L$ is a {\it maximal} prefix code. 

\bigskip

\noindent 
{\bf Theorem 4.5B (Equivalent} \textbf{\textsl{normal}} {\bf morphism).} 
{\it For every right-ideal morphism $\varphi$ with finite domain code
there exists a {\rm normal} right-ideal morphism $\varphi_0$ with finite 
domain code, such that $\varphi = \varphi_0$ in $M_{k,1}$. 
Moreover, 

\smallskip

$|{\rm imC}(\varphi_0)| \ = \ |\varphi_0({\rm domC}(\varphi_0))|$ 
$\ \le \ $ 
$|A| \cdot (\ell(\varphi(P)) + 1) \cdot |\varphi(P)|$,

\smallskip

$|{\rm domC}(\varphi_0)|$ $\ \le \ $
$|P| \cdot |A| \cdot (\ell(\varphi(P)) + 1) \cdot |\varphi(P)|$,

\smallskip

$\ell({\rm imC}(\varphi_0)) \, = \, \ell(\varphi_0({\rm domC}(\varphi_0)))$
$\, = \, \ell(\varphi({\rm domC}(\varphi)))$,

\smallskip

$\ell({\rm domC}(\varphi_0)) \ \le \ $
$\ell({\rm domC}(\varphi)) \, + \, \ell(\varphi({\rm domC}(\varphi)))$.  
}
 
\bigskip

\noindent {\bf Proof.} Let $P = {\rm domC}(\varphi)$, 
$Q = {\rm imC}(\varphi)$, $P_0 = {\rm domC}(\varphi_0)$, 
$Q_0 = {\rm imC}(\varphi_0)$. 
For each $p \in P$, let $\varphi(p) \, W_{\varphi(p)}$ be the the set of 
leaves of the saturated tree 
$\, {\rm s}T\big(\varphi(p), \, \varphi(P) \cap \varphi(p) \, A^*\big)$.
By Remark (3) above, $W_{\varphi(p)}$ is a finite maximal prefix code.   
Now we define $\varphi_0$ as follows:

\smallskip

 \ \ \  \ \ \   $\varphi_0$ is the restriction of $\varphi$ to  
$ \ \bigcup_{p \in P} \, p \, W_{\varphi(p)} \, A^*$. 

\smallskip

\noindent Let us verify that $\varphi_0$ has the required properties.
Since $P$ and $W_{\varphi(p)}$ are finite prefix codes, 
$\bigcup_{p \in P} p \, W_{\varphi(p)}$ is a finite prefix code. So,

\smallskip

 \ \ \  \ \ \   
${\rm domC}(\varphi_0) \ = \ \bigcup_{p \in P} \, p \, W_{\varphi(p)}$. 

\smallskip

\noindent
Since each $W_{\varphi(p)}$ is a maximal prefix code, the right ideal 
$\, \bigcup_{p \in P} p \, W_{\varphi(p)} \, A^*$ is essential in the right 
ideal $P A^*$; hence $\varphi$ and $\varphi_0$ are equal as elements of 
$M_{k,1}$.
Finally, let us show that $\varphi_0({\rm domC}(\varphi_0))$ is a prefix 
code.  We have 

\smallskip

 \ \ \  \ \ \ $\varphi_0({\rm domC}(\varphi_0))$ $ \ = \ $ 
$\bigcup_{p \in P} \, \varphi(p) \ W_{\varphi(p)}$, 

\smallskip

\noindent which is the set of leaves of the union of the saturated prefix 
trees 
$\, {\rm s}T\big(\varphi(p), \, \varphi(P) \cap \varphi(p) \, A^*\big)$, 
for $p$ ranging over $P$.  For each $p \in P$, the leaves of 
$\, {\rm s}T\big(\varphi(p), \, \varphi(P) \cap \varphi(p) \, A^*\big)$ form
the prefix code $\varphi(p) \, W_{\varphi(p)}$. 
For $p_1 \neq p_2$ in $P$, if $\varphi(p_1)$ is a prefix of $\varphi(p_2)$ 
then the leaves of 
$\, {\rm s}T\big(\varphi(p_2), \, \varphi(P) \cap \varphi(p_2) \, A^*\big)$ 
are a subset of the leaves of 
$\, {\rm s}T\big(\varphi(p_1), \, \varphi(P) \cap \varphi(p_1) \, A^*\big)$,
so the union of these two leaf sets is just the leaf set of 
$\, {\rm s}T\big(\varphi(p_1), \, \varphi(P) \cap \varphi(p_1) \, A^*\big)$;
a similar thing happens if $\varphi(p_2)$ is a prefix of $\varphi(p_1)$.
So in $\, \bigcup_{p \in P} \, \varphi(p) \ W_{\varphi(p)} \, $ we can 
ignore elements $p$ of $P$ for which $\varphi(p)$ is a strict prefix of 
another element of $\varphi(P)$.
If $\varphi(p_1)$ and $\varphi(p_2)$ are not prefix-comparable,
then the leaves of 
$\, {\rm s}T\big(\varphi(p_i), \, \varphi(P) \cap \varphi(p_i) \, A^*\big)$ 
have $\varphi(p_i)$ as a prefix, so these two trees have leaf sets that are
two-by-two prefix-incomparable (namely the sets 
$\varphi(p_1) \, W_{\varphi(p_1)}$ and $\varphi(p_2) \, W_{\varphi(p_2)}$). 
The union of prefix codes that are two-by-two prefix-incomparable forms a 
prefix code; hence,
$\, \bigcup_{p \in P} \, \varphi(p) \ W_{\varphi(p)} \, $ is a prefix code. 

Now, since $\varphi_0({\rm domC}(\varphi_0))$ is a prefix code it follows
that $\, {\rm imC}(\varphi_0) = \varphi_0({\rm domC}(\varphi_0))$, so  
$\varphi_0$ is normal.  This proves the first part of the theorem.

\smallskip

Let us prove the formulas.
We saw that $\, {\rm imC}(\varphi_0) = \varphi_0({\rm domC}(\varphi_0))$ 
$=$ $\, \bigcup_{p \in P} \, \varphi(p) \ W_{\varphi(p)}$, and 
$\varphi(p) \ W_{\varphi(p)}$ is the leaf set of the saturated tree 
$\, {\rm s}T\big(\varphi(p), \varphi(P) \cap \varphi(p) A^*\big)$.
By the definition of prefix trees, the vertices of all the (non-saturated) 
trees $\, T\big(\varphi(p), \varphi(P) \cap \varphi(p) A^*\big) \,$ are 
subsets of ${\sf pref}(\varphi(P))$. By Remark (2) above, the number of
leaves in a saturated tree
$\, {\rm s}T\big(\varphi(p), \varphi(P) \cap \varphi(p) A^*\big) \,$ 
is at most $|A|$ times the number of vertices of the non-saturated tree.
Hence, $|{\rm imC}(\varphi_0)| \le |A| \cdot |{\sf pref}(\varphi(P))|$.
Moreover, for any finite $Z \subset A^*$,  
$\, |{\sf pref}(Z)| \le (1 + \ell(Z)) \cdot |Z|$, hence,
$|{\rm imC}(\varphi_0)| \le$
$|A| \cdot (\ell(\varphi(P)) + 1) \cdot |\varphi(P)|$.

We have $\, {\rm domC}(\varphi_0)$ $=$
$\bigcup_{p \in P} p \, W_{\varphi(p)}$, and
$\varphi(p) \, W_{\varphi(p)}$ is the leaf set of 
$\, {\rm s}T\big(\varphi(p), \varphi(P) \cap \varphi(p) A^*\big)$. 
Hence by the same reasoning as for $|{\rm imC}(\varphi_0)|$: 
$\, |W_{\varphi(p)}| = |\varphi(p) \, W_{\varphi(p)}|$ $\, \le \,$
$|A| \cdot (\ell(\varphi(P)) + 1) \cdot |\varphi(P)|$.  Hence, 
$|{\rm domC}(\varphi_0)|$ $\, \le \,$ $\sum_{p \in P} |W_{\varphi(p)}|$ 
$\, \le \,$
$\sum_{p \in P} \, |A| \cdot (\ell(\varphi(P)) + 1) \cdot |\varphi(P)|$
$\, \le \, $ 
$|P| \cdot |A| \cdot (\ell(\varphi(P)) + 1) \cdot |\varphi(P)|$. 

We have 
$\, {\rm imC}(\varphi_0) = \bigcup_{p \in P} \varphi(p) \, W_{\varphi(p)}$,
and $\varphi(p) \, W_{\varphi(p)}$ is the leaf set of 
$\, {\rm s}T\big(\varphi(p), \, \varphi(P) \cap \varphi(p) \, A^*\big)$.
Hence, $\ell({\rm imC}(\varphi_0)) \le \ell(\varphi(P))$; indeed,
tree saturation does not increase the depth of a tree, and the depth of
$\, T\big(\varphi(p), \, \varphi(P) \cap \varphi(p) \, A^*\big) \, $ is
$\, \le \ell(\varphi(P))$.

We have $\, {\rm domC}(\varphi_0) = \bigcup_{p \in P} p \, W_{\varphi(p)}$.
And $\ell(W_{\varphi(p)}) \le \ell(\varphi(p) \, W_{\varphi(p)})$ $\le$ $\ell(\varphi(P))$, since 
$\varphi(p) \, W_{\varphi(p)}$ is the leaf set of
$\, {\rm s}T\big(\varphi(p), \, \varphi(P) \cap \varphi(p) \, A^*\big)$.
Hence, for every $x \in {\rm domC}(\varphi_0)$ we have 
$x \in p W_{\varphi(p)}$ for some $p \in P$, so 
$\, |x| \le |p| + \ell(W_{\varphi(p)})$. Therefore, 
$\ell({\rm domC}(\varphi_0)) \, \le \, \ell(P) + \ell(\varphi(P))$. 
  \ \ \ $\Box$

\bigskip

\noindent Theorem 4.5B tells us that as far as $M_{k,1}$ is concerned, all 
right-ideal morphisms are normal. 
\footnote{\, The concept of normal morphism and Theorem 4.5 enable 
us to rehabilitate the {\it image code formula} (which was incorrect as
stated in Theorem 4.5 of \cite{BiMonTh}, but which is correct when one adds
the hypothesis that the morphisms $\varphi_i$ are {\em normal}).} 

\bigskip

\noindent {\bf Corollary 4.5C (Image code formula).} \\   
{\it Let $\varphi_i$ be a right-ideal morphism (for $i = 1, \ldots, n$), and 
let $\, \Phi =  \varphi_n \circ \ \ldots \ \circ \varphi_1$. 

\medskip

\noindent {\bf (1)} If $\varphi_i$ is {\rm normal} for $2 \le i \le n$, 
then

\medskip

\hspace{.3in}   $|{\rm imC}(\Phi)|$ $ \ \leq \ $
$\sum_{i=1}^n |{\rm imC}(\varphi_i)|$ .

\medskip

\noindent {\bf (2)} If all $\varphi_i$ are {\rm normal} (for
$1 \le i \le n$), then

\medskip

\hspace{.3in}    $\ell({\rm domC}(\Phi) \cup {\rm imC}(\Phi)) \ \leq \ $
  $\sum_{i=1}^n \ell({\rm domC}(\varphi_i) \cup {\rm imC}(\varphi_i))$ .
}

\bigskip

\noindent {\bf Proof.} {\bf (1)} follows immediately from Theorem 
\ref{sum_of_imC}(1), and {\bf (2)} follows from \ref{sum_of_imC}(2) and 
\ref{sum_of_imC}(4).  
  \ \ \ $\Box$

\bigskip

\noindent {\bf Counter-examples:} 

\smallskip

\noindent {\bf (1)}
The following shows that the image code formula of Corollary 4.5C(1) 
is wrong in some examples when $\varphi_2$ is not normal (but $\varphi_1$
is normal).  Let $A = \{0,1\}$, $n \ge 2$, and

\smallskip

$\varphi_1 \ = \ \{(01, 00), \, (00, 01), \, (10, 1011), \, (11, 1100)\}$, 
  \ and

\smallskip

$\varphi_2 \ = \ \{(00u0, 000u1) : u \in \{0,1\}^{n-1}\}$ $\ \cup \ $
  $\{(01v0, 001v1) : v \in \{0,1\}^{n-1}\}$ $\ \cup \ $ 
  $\{(10, 000), \, (11, 001)\}$.

\medskip

\noindent So, ${\rm imC}(\varphi_1) = \{00, 01, 1011, 1100\}$, and 
${\rm imC}(\varphi_2) = \{000, 001\}$, hence 
$|{\rm imC}(\varphi_1)| + |{\rm imC}(\varphi_2)| = 6$.
Note that the right-ideal morphisms $\varphi_1$ and $\varphi_2$ are in 
maximally extended form.

Now, $\varphi_2 \circ \varphi_1: \ 01u0 \mapsto 00 u0 \mapsto 000 u1$
and $\varphi_2 \circ \varphi_1: \ 00 v0 \mapsto 01 v0 \mapsto 001 v1$, 
for all $u, v \in \{0,1\}^{n-1}$;  and
$\varphi_2 \circ \varphi_1: \ 10 \mapsto 1011 \mapsto  00011$, 
$\varphi_2 \circ \varphi_1: \ 11 \mapsto 1100 \mapsto 00100$.
Note that $\varphi_2 \circ \varphi_1$ is in maximally extended form. 

Then $\, {\rm imC}(\varphi_2 \circ \varphi_1) = \{00011, 00100\}$ $\cup$
$000 \, \{00, 01,11\} \, \{0,1\}^{n-2}$ $\cup$ 
$001 \, \{00, 01,11\} \, \{0,1\}^{n-2}$. Thus when $n \ge 2$:  
 \ $2 + 6 \cdot 2^{n-2} = |{\rm imC}(\varphi_2 \circ \varphi_1)|$
$ \ \not\le \ $ $|{\rm imC}(\varphi_1)| + |{\rm imC}(\varphi_2)| = 6$. 
 \ \ \ $\Box$

\medskip

\noindent {\bf (2)} The following shows that the formula of Corollary 
4.5C(2) is wrong in some examples when $\varphi_2$ is not normal (but 
$\varphi_1$ is normal). We abbreviate 
$\ell({\rm domC}(\varphi) \cup {\rm imC}(\varphi))$ by $\ell(\varphi)$.  
Let $A = \{0,1\}$, $n \ge 2$, and

\smallskip

$\varphi_1 \ = \ \{(0, 0^n)\}$, \ and

\smallskip

$\varphi_2 \ = \ \{(0, 0^{n+1}), \, (1, 0)\}$.

\medskip

\noindent So, $\ell(\varphi_1) = n$, and $\ell(\varphi_2) = 1$ since 
${\rm imC}(\varphi_2) = \{0\}$.  Now, 
$\varphi_2 \circ \varphi_1 = \{(0, 0^{2n})\}$. Thus when
$n \ge 2$: \ $2n = \ell(\varphi_2 \circ \varphi_1) \ \not\le \ $
$\ell(\varphi_2) + \ell(\varphi_1) = n + 1$.
 \ \ \ $\Box$

\bigskip

For elements of ${\it Inv}_{k,1}$ the image code has the same size as the 
domain code, which is also the table size. Moreover, injective right-ideal
morphisms are normal, thus Corollary 4.5C implies:

\begin{cor} \     
For all injective right-ideal morphisms $\varphi, \psi$:  
 \ \ $\|\psi \circ \varphi\| \leq \|\psi\| + \|\varphi\|$.
 \ \ \ \ \ $\Box$
\end{cor}
In other words, the table size formula holds for ${\it Inv}_{k,1}$.
Another immediate consequence of Theorem \ref{sum_of_imC} is the following.

\begin{cor} \label{repeated_sum_of_imC} \   
Let $\varphi_i$ be {\em normal} right-ideal morphisms for 
$i = 1, \ldots, n$, and let $c_1, c_2$ be positive constants.

\medskip

\noindent {\bf (1)} \ \ \ \ If $\, |{\rm imC}(\varphi_i)| \leq c_1$ for 
all $i$ then \    
$|{\rm imC}(\varphi_n \circ \ldots \circ \varphi_1)| \ \le \ c_1 \, n$.

\medskip

\noindent {\bf (2)} \ \ \ \ If $\, \ell({\rm imC}(\varphi_i)) \le c_2$ 
for all $i$ then \    
$\ell({\rm imC}(\varphi_n \circ \ldots \circ \varphi_1)) \ \le \ c_2 \, n$.
  \ \ \  \ \ \  \ $\Box$
\end{cor}

The position transposition $\tau_{i,j}$ (with $0<i<j$) is, by definition, 
the partial permutation of $A^*$ which transposes the letters at positions 
$i$ and $j$; $\tau_{i,j}$ is undefined on words of length $< j$. More
precisely, we have \
${\rm domC}(\tau_{i,j}) = {\rm imC}(\tau_{i,j}) = A^j$, and 
 \ $u \alpha v \beta \mapsto u \beta v \alpha$ \ for all letters
$\alpha, \beta \in A$ and all words $u \in A^{i-1}$ and $v \in A^{j-i-1}$.
In this form, $\tau_{i,j}$ is equal to its maximum essential extension.

\begin{cor} \label{wordlength_tau} \
The word-length of $\tau_{i,j}$ over any finite generating set of $M_{k,1}$
is exponential.
\end{cor}
{\bf Proof.} We have $|{\rm imC}(\tau_{i,j})| = k^j$. The Corollary follows
then from Corollary \ref{repeated_sum_of_imC}(1).
 \ \ \ $\Box$


\subsection{Some algorithmic problems about right-ideal morphisms}

We consider several problems about right-ideal morphisms of $A^*$
and show that they have deterministic polynomial-time algorithms. We also
show that the word problem of $M_{k,1}$ over $\Gamma_{k,1} \cup $
$\{ \tau_{i,i+1} : 0 < i \}$ is {\sf coNP}-complete, where $\Gamma_{k,1}$ 
is any finite generating set of $M_{k,1}$.  We saw that $\Gamma_{k,1}$ can
be chosen so as to consist of normal right-ideal morphisms.

\begin{lem} \label{essent_intersect_algo} \
There are deterministic polynomial time algorithms for the following
problems. \\
{\rm Input:} Two finite prefix codes $P_1, P_2 \subset A^*$, given explicitly
by lists of words. \\
{\rm Output 1:} The finite prefix code $\Pi \subset A^*$ such that
$\Pi A^* = P_1A^* \cap P_2A^*$, where $\Pi$ is produced explicitly
as a list of words. \\
{\rm Question 2:} Is $P_1A^* \cap P_2A^*$ essential in $P_1A^*$ (or in
$P_2A^*$, or in both)?
\end{lem}
{\bf Proof.}  We saw already that $\Pi$ exists and
$\Pi \subseteq P_1 \cup P_2$; see Lemma 3.3 of \cite{BiThomps} (quoted before
Lemma \ref{Im_of_rightideal} above).

\smallskip

\noindent {\sf Algorithm for Output 1:}\  Since 
$\Pi \subseteq P_1 \cup P_2$, we just need to search for
the elements of $\Pi$ within $P_1 \cup P_2$.
For each $x \in P_1$ we check whether $x$ also
belongs to $P_2A^*$ (by checking whether any element of $P_2$ is a prefix of
$x$). Since $P_1$ and $P_2$ are explicitly given as lists, this takes
polynomial time. Similarly, for each $x \in P_2$ we check whether $x$ also
belongs to $P_1A^*$. Thus, we have computed the set \,
$\Pi_1 = (P_1 \cap P_2A^*) \, \cup \, (P_2 \cap P_1A^*)$.
Now, $\Pi$ is obtained from $\Pi_1$ by eliminating every word that has
another word of $\Pi_1$ as a prefix. Since $\Pi_1$ is explicitly listed,
this takes just polynomial time.

\smallskip

\noindent {\sf Algorithm for Question 2:} \ We first compute $\Pi$ by the
previous algorithm. Next, we check whether every $p_1 \in P_1$ is a prefix 
of some $r \in \Pi$; since $P_1$ and $\Pi$ are given by explicit lists, 
this takes just polynomial time. For $P_2$ it is similar.
 \ \ \ $\Box$

\begin{lem} \label{varphi_of_set} \
The following input-output problem has a deterministic polynomial-time
algorithm. \\
$\bullet$ {\rm Input:} A finite set $S \subset A^*$, and $m$ right-ideal 
morphisms $\psi_j$ for $j = 1, \ldots, m$, where $S$ is given by an
explicit list of words, and each $\psi_j$ is given explicitly by the list of 
pairs of words $\{(x, \psi_j(x)) : x \in {\rm domC}(\psi_j)\}$.  \\
$\bullet$ {\rm Output:} The finite set $\, \psi_m \ldots \psi_1(S)$, given 
explicitly by a list of words.
\end{lem}
{\bf Proof.} Let $\Psi = \, $ 
$\psi_m \circ \ \ldots \ \circ \psi_1 \circ {\sf id}_S$.
Then $\psi_m \ldots \psi_1(S) = \Psi({\rm domC}(\Psi))$.
By Theorem \ref{sum_of_imC}(3) and (5), $\ell(\Psi({\rm domC}(\Psi)))$ 
$\le$ $\ell(S) + \sum_{i=1}^m \ell(\psi_i({\rm domC}(\psi_i))) \, $ and 
$\, |\Psi({\rm domC}(\Psi))| \le$
$|S| + \sum_{i=1}^m |\psi_i({\rm domC}(\psi_i))|$.
So the size of $\psi_m \ldots \psi_1(S)$, in terms of the number of words 
and their lengths, is polynomially bounded by the size of the input.

We now compute $\psi_m \ldots \psi_1(S)$ by applying $\psi_j$ to 
$\psi_{j-1} \ldots \psi_1(S)$ for increasing $j$. Since the sizes of the 
sets remain polynomially bounded, this algorithm takes polynomial time. 
   \ \ \ $\Box$

\begin{cor} \label{computing_imC} \
The following input-output problems have deterministic polynomial-time
algorithms. \\
$\bullet$ {\rm Input:} A list of $n$ right-ideal morphisms
$\varphi_i$ for $i = 1, \ldots, n$, given explicitly by finite tables.  \\
$\bullet$ {\rm Output 1:} A finite set, as an explicit list of words, that 
contains $\, \varphi_n \ldots $
$\varphi_1({\rm domC}(\varphi_n \ldots \varphi_1))$. \\   
$\bullet$ {\rm Output 2:} The finite set 
$\, {\rm imC}(\varphi_n \ldots \varphi_1)$, as an explicit list of words.  
\end{cor}
{\bf Proof.} (1) By Theorem \ref{sum_of_imC}(5) we have
$\, \varphi_n \ldots \varphi_1({\rm domC}(\varphi_n \ldots \varphi_1))$
$ \ \subseteq \ \bigcup_{i=1}^n$
$\varphi_n \ldots \varphi_i({\rm domC}(\varphi_i))$. By Lemma 
\ref{varphi_of_set}, each set 
$\varphi_n \ldots \varphi_i({\rm domC}(\varphi_i))$, as well as their union,
is computable in polynomial time (as an explicit list of words).  

(2) Let $\Phi = \varphi_n \ldots \varphi_1$, 
$P_i = {\rm domC}(\varphi_i)$, and $Q_i = {\rm imC}(\varphi_i)$.
As in the proof of Theorem \ref{sum_of_imC}(1), 
 \ ${\rm Dom}(\varphi_2 \circ \varphi_1) = $ 
$\varphi_1^{-1}(Q_1A^* \cap P_2A^*)$, 
 \ ${\rm Im}(\varphi_2 \circ \varphi_1) = \varphi_2(Q_1A^* \cap P_2A^*)$, 
and the maps 
$\varphi_1^{-1}(Q_1A^* \cap P_2A^*) \ \stackrel{\varphi_1}{\longrightarrow}$
$ \ Q_1A^* \cap P_2A^* \stackrel{\varphi_2}{\longrightarrow} \ $
$\varphi_2(Q_1A^* \cap P_2A^*)$ \ are total and onto.
By Lemma 3.3 of \cite{BiThomps} (mentioned before Theorem \ref{sum_of_imC}) 
we have $Q_1A^* \cap P_2A^* = S_1 A^*$ for some finite prefix code $S_1$ 
with $S_1 \subseteq Q_1 \cup P_2$. Moreover, by Lemma \ref{Im_of_rightideal}, 
$\varphi_2(S_1 A^*) = R_2 A^*$, where 
$\, {\rm imC}(\varphi_2 \varphi_1) = R_2 \subseteq \varphi_2(S_1)$.  

By induction, for $j \ge 2$ suppose 
$\, {\rm imC}(\varphi_j \ldots \varphi_1) = R_j$ $\subseteq$ 
$\varphi_j(S_{j-1})$, where $R_j$ and $S_{j-1}$ are finite prefix codes 
such that $S_{j-1} \subseteq R_{j-1} \cup P_j$, 
 \ $S_{j-1} A^* = R_{j-1} A^* \cap P_j A^*$, \ $R_j A^* = $
${\rm Im}(\varphi_j \ldots \varphi_1) = \varphi_j(S_{j-1} A^*)$, 
and the maps $\varphi_j^{-1}(R_j A^* \cap P_{j+1} A^*)$ 
$ \ \stackrel{\varphi_j}{\longrightarrow}$
$ \ R_j A^* \cap P_{j+1} A^* \stackrel{\varphi_{j+1}}{\longrightarrow} \ $
$\varphi_{j+1}(R_j A^* \cap P_{j+1} A^*)$ \ are total and onto.
Then by Lemma 3.3 of \cite{BiThomps} we again have 
$R_j A^* \cap P_{j+1} A^* = S_j A^*$ for some finite prefix code $S_j$ 
with $S_j \subseteq R_j \cup P_{j+1}$; and by Lemma \ref{Im_of_rightideal},
$\varphi_{j+1}(S_j A^*) = R_{j+1} A^*$ for some finite prefix code 
$R_{j+1}$ such that ${\rm imC}(\varphi_{j+1} \varphi_j \ldots \varphi_1) = $
$R_{j+1} \subseteq \varphi_{j+1}(S_j)$.
Applying Theorem \ref{sum_of_imC} to 
$\, R_i = {\rm imC}(\varphi_i \ \ldots \ \varphi_1)$ for any $i \ge 2$
we have 

\smallskip

$\, |R_i| \, \le \, |\varphi_i(P_i)| + \ \ldots \ $ $+$
$|\varphi_2(P_2)| + |{\rm imC}(\varphi_1)|$, \ and

\smallskip

$\, \ell(R_i) \, \le \, \ell(\varphi_i(P_i)) + \ \ldots \ $ $+$
$\ell(\varphi_2(P_2)) + \ell({\rm imC}(\varphi_1))$. 

\smallskip

\noindent Since $S_j \subseteq P_j \cup R_{j-1}$, we have 
$|S_j| \le |P_j| + |R_{j-1}|$ $\, \le \,$ 
$|P_j| + |\varphi_{j-1}(P_{j-1})| + \ \ldots \ $ $+$
$|\varphi_2(P_2)| + |{\rm imC}(\varphi_1)|$, \ and
$\ell(S_j) \le \ell(P_j) + \ell(R_{j-1})$ $\, \le \,$
$\ell(P_j) + \ell(\varphi_{j-1}(P_{j-1})) + \ \ldots \ $ $+$ 
$\ell(\varphi_2(P_2)) + \ell({\rm imC}(\varphi_1))$.
Thus, the size of each $R_i$ and $S_j$ is less than the input size;  
by input size we mean the total length of all the words in the input 
lists.

By Lemma \ref{essent_intersect_algo}, the prefix code $S_j$ is computed from
$R_j$ and $P_{j+1}$, as an explicit list, in time 
$\le T_j(|P_j| + \ell(P_j) + |R_{j-1}| + \ell(R_{j-1}))$, for some 
polynomial $T_j(.)$.
And $R_{j+1}$ is computed from $S_j$ by applying $\varphi_{j+1}$ to $S_j$, 
and then keeping the elements that do not have a prefix in 
$\varphi_{j+1}(S_j)$. Computing $\varphi_{j+1}(S_j)$ takes at most
quadratic time, and finding the prefix code in $\varphi_{j+1}(S_j)$ also
takes at most quadratic time. 

In the end we obtain $R_n = {\rm imC}(\varphi_n \ldots \varphi_1)$ 
as an explicit list of words.  
 \ \ \ $\Box$

\bigskip

When we consider the word problem of $M_{k,1}$ over a finite generating
set, we measure the input size by the length of input word (with each
generator having length 1).
But for the word problem of $M_{k,1}$ over the infinite generating 
set $\Gamma_{k,1} \cup \{\tau_{i-1,i} : i > 1\}$ we count the length of
the position transpositions $\tau_{i-1,i}$ as $i$,
in the definition of the input size of the word problem. 
Indeed, at least $\log_2 i$ bits are needed to describe the subscript $i$
of $\tau_{i-1,i}$. Moreover, in the connection between $M_{k,1}$ (over
$\Gamma_{k,1} \cup \{\tau_{i-1,i} : i > 1\}$) and circuits, $\tau_{i-1,i}$ 
is interpreted as the wire-crossing operation of wire number $i$ and wire
number $i-1$; this suggests that viewing the size of $\tau_{i-1,i}$ as $i$ 
is more natural than $\log_2 i$. In any case, we will see next that the word 
problem of $M_{k,1}$ over $\Gamma_{k,1} \cup \{\tau_{i-1,i} : i > 1\}$ is 
{\sf coNP}-complete, even if the size of $\tau_{i-1,i}$ is more generously
measured as $i$; this is a stronger result than if $\log_2 i$ were used.

\begin{thm} \label{wp_in_coNP} {\bf ({\sf coNP}-complete word problem).} \    
The word problem of $M_{k,1}$ over the infinite generating set 
$\Gamma_{k,1} \cup \{\tau_{i-1,i} : i > 1\}$ is {\sf coNP}-complete, 
where $\Gamma_{k,1}$ is any finite generating set of $M_{k,1}$.  
\end{thm}
{\bf Proof.}  In \cite{BiCoNP} (see also \cite{BiFact}) it was shown that 
the word problem of the Thompson-Higman group $G_{k,1}$ over 
$\Gamma_{G_{k,1}} \cup \{ \tau_{i-1,i} : i > 1\}$ \ is {\sf coNP}-complete,
where $\Gamma_{G_{k,1}}$ is any finite generating set of $G_{k,1}$. Hence, 
since the elements of the finite set $\Gamma_{G_{k,1}}$ can be expressed by
a finite set of words over $\Gamma_{k,1}$, it follows that  the 
word problem of $M_{k,1}$ over $\Gamma_{k,1} \cup \{\tau_{i-1,i} : i > 1\}$ 
is {\sf coNP}-hard.

We will prove now that the word problem of $M_{k,1}$ over 
$\Gamma_{k,1} \cup \{\tau_{i-1,i} : i > 1\}$ belongs to {\sf coNP}. 
The {\it input}
of the problem consists of two words $(\rho_m, \ldots, \rho_1)$ and 
$(\sigma_n, \ldots, \sigma_1)$ over 
$\Gamma_{k,1} \cup \{\tau_{i-1,i} : i > 1\}$.  The {\it input size} is the 
weighted length of the words $(\rho_m, \ldots, \rho_1)$ and
$(\sigma_n, \ldots, \sigma_1)$, where each generator in $\Gamma_{k,1}$ has 
weight 1, and each generator of the form $\tau_{i-1,i}$ has weight $i$.
For every right-ideal morphism $\varphi$ we abbreviate 
$\ell({\rm domC}(\varphi) \cup \varphi({\rm domC}(\varphi)))$ by 
$\ell(\varphi)$; recall that for a finite set $X \subset A^*$, $\ell(X)$
denotes the length of a longest word in $X$.

Since $\Gamma_{k,1}$ is finite there is a constant $c>0$ such that 
$c \geq \ell(\gamma)$ for all $\gamma \in \Gamma_{k,1}$; also, for each 
$\tau_{i-1,i}$ we have $\ell(\tau_{i-1,i}) = i$. 
By Theorem \ref{sum_of_imC}, the table of 
$\sigma_n \circ \ldots \circ \sigma_1$ (and more generally, the table of
$\sigma_j \circ \ldots \circ \sigma_1$ for any $j$ with $n \geq j \geq 1$)
contains only words of length $\leq \sum_{j=1}^n \ell(\sigma_j)$, and 
similarly for $\rho_m \circ \ldots \circ \rho_1$ (and for 
$\rho_i \circ \ldots \circ \rho_1$, $m \geq i \geq 1$). 
So all the words in the tables for any 
$\sigma_j \circ \ldots \circ \sigma_1$  and any
$\rho_i \circ \ldots \circ \rho_1$ have lengths that are linearly bounded 
by the size of the input 
$\big((\rho_m, \ldots, \rho_1)$, $(\sigma_n, \ldots, \sigma_1)\big)$. 

\medskip

\noindent {\sf Claim.} {\it Let $N = \max\{\sum_{i=1}^m \ell(\rho_i), \ $
  $\sum_{j=1}^n \ell(\sigma_j)\}$. Then \ 
$\rho_m \cdot \ldots \cdot \rho_1 \neq \sigma_n \cdot \ldots \cdot \sigma_1$
 \ in $M_{k,1}$ \ iff \ there exists $x \in A^N$ such that 
 \ $\rho_m \circ \ldots \circ \rho_1(x)$ $\neq$ 
$\sigma_n \circ \ldots \circ \sigma_1(x)$.          }
 
\smallskip

\noindent {\sf Proof of the Claim:} As we saw above, the tables of  
$\rho_m \circ \ldots \circ \rho_1$ and $\sigma_n \circ \ldots \circ \sigma_1$
only contain words of length $\leq N$. Thus, restricting 
$\rho_m \circ \ldots \circ \rho_1$ and $\sigma_n \circ \ldots \circ \sigma_1$
to $A^N A^*$ is an essential restriction, and the resulting tables have
domain codes in $A^N$. Therefore, $\rho_m \cdot \ldots \cdot \rho_1$
and $\sigma_n \cdot \ldots \cdot \sigma_1$ are equal (as elements of
$M_{k,1}$) \, iff \,  
$\rho_m \circ \ldots \circ \rho_1$ and $\sigma_n \circ \ldots \circ \sigma_1$ 
are equal on $A^N$. \ \ \ \ \ [End, Proof of Claim] 

\medskip

The number $N$ in the Claim is immediately obtained form the input. 
Based on the Claim, we obtain a nondeterministic polynomial-time algorithm 
which decides (nondeterministically) whether there exists $x \in A^N$ such 
that \ $\rho_m \circ \ldots \circ \rho_1(x)$
$ \neq \sigma_n \circ \ldots \circ \sigma_1(x)$, as follows:

The algorithm guesses $x \in A^N$, computes 
$\rho_m \circ \ldots \circ \rho_1(x)$ and
$\sigma_n \circ \ldots \circ \sigma_1(x)$, and checks that they are different
words ($\in A^*$) or that one is undefined and the other is a word. 
Applying Theorem \ref{sum_of_imC} to 
$\, \rho_m \circ \ldots \circ \rho_1 \circ {\sf id}_{A^N} \, $ and to 
$\, \sigma_n \circ \ldots \circ \sigma_1 \circ {\sf id}_{A^N} \, $ shows 
that  \ $|\rho_m \circ \ldots \circ \rho_1(x)| \leq 2N$ \ and 
 \ $|\sigma_n \circ \ldots \circ \sigma_1(x)| \leq 2N$; here 
$|\rho_m \circ \ldots \circ \rho_1(x)|$ denotes the length of the word
$\rho_m \circ \ldots \circ \rho_1(x) \in A^*$, and similarly for 
$\sigma_n \circ \ldots \circ \sigma_1(x)$.
Also by Theorem \ref{sum_of_imC}, all intermediate results 
(as we successively apply $\rho_i$ for $i = 1, \ldots, m$, or $\sigma_j$
for $j = 1, \ldots, n$) are words of length $\leq 2N$. 
These successive words are computed by applying the table of $\rho_i$ or
$\sigma_j$ (when $\rho_i$ or $\sigma_j$ belong to $\Gamma_{k,1}$), or by
directly applying the position permutation $\tau_{h,h-1}$ (if $\rho_i$ or
$\sigma_j$ is $\tau_{h,h-1}$).  
Thus, the output $\rho_m \circ \ldots \circ \rho_1(x)$ (and similarly, 
$\sigma_n \circ \ldots \circ \sigma_1(x)$) can be computed in polynomial
time.
 \ \ \ $\Box$


\subsection{The word problem of $M_{k,1}$ is in {\sf P}}

We now move ahead with the the proof of our main result.

\begin{thm} \label{wp_M_in_P} {\bf (Word problem in {\sf P}).} \  
The word problem of the Thompson-Higman monoids $M_{k,1}$, over any finite 
generating set, can be decided in deterministic polynomial time.
\end{thm}

We assume that a fixed finite generating set $\Gamma_{k,1}$ of $M_{k,1}$ has
been chosen. The input consists of two sequences $(\rho_m, \ldots, \rho_1)$ 
and $(\sigma_n, \ldots, \sigma_1)$ over $\Gamma_{k,1}$, and the input size
is $m+n$; since $\Gamma_{k,1}$ is finite and fixed, it does not matter
whether we choose $m+n$ as input size, or the sum of the lengths of all the
words in the tables of the elements of $\Gamma_{k,1}$.
We want to decide in deterministic polynomial time whether, as elements of
$M_{k,1}$, the products $\rho_m \cdot \ldots \cdot \rho_1$ and 
$\sigma_n \cdot \ldots \cdot \sigma_1$ are equal.

\medskip

\noindent {\bf Overview of the proof:}

\smallskip

\noindent $\bullet$ \ We compute the finite sets
$\, {\rm imC}(\rho_m \circ \ldots \circ \rho_1)$, 
$\, {\rm imC}(\sigma_n \circ \ldots \circ \sigma_1) \subset A^*$, explicitly 
described by lists of words. By Corollary \ref{computing_imC} (Output 2)
this can be done in polynomial time, and these sets have polynomial size. 
(Note however that by Proposition \ref{exp_table}, the table sizes of
$\rho_m \circ \ldots \circ \rho_1$ or $\sigma_n \circ \ldots \circ \sigma_1$
could be exponential in $m$ or $n$.) 

\smallskip

\noindent $\bullet$ \ We check whether 
 ${\rm Im}(\rho_m \circ \ldots \circ \rho_1) $ $\cap $
 ${\rm Im}(\sigma_n \circ \ldots \circ \sigma_1)$ is essential in 
${\rm Im}(\rho_m \circ \ldots \circ \rho_1)$ and in
${\rm Im}(\sigma_n \circ \ldots \circ \sigma_1)$. By Lemma 
\ref{essent_intersect_algo} (Question 2) this can be done in polynomial 
time.  If the answer is ``no'' then $\rho_m \cdot \ldots \cdot \rho_1$ 
$\neq$ $\sigma_n \cdot \ldots \cdot \sigma_1$ in $M_{k,1}$, since they 
don't have a common maximum essential extension. Otherwise, the computation
continues.

\smallskip

\noindent $\bullet$ \ We compute the finite prefix code $\Pi \subset A^*$ 
such that $\Pi A^* = {\rm Im}(\rho_m \circ \ldots \circ \rho_1) $ $\cap $
 ${\rm Im}(\sigma_n \circ \ldots \circ \sigma_1)$. By Lemma
\ref{essent_intersect_algo} (Output 1) this can be done in polynomial time, 
and $\Pi$ has polynomial size. Hence, the table of ${\sf id}_{\Pi A^*}$ 
can be computed in polynomial time.

\smallskip

\noindent $\bullet$ \ We restrict $\rho_m \circ \ldots \circ \rho_1$ and 
$\sigma_n \circ \ldots \circ \sigma_1$ in such a way that their images are 
in $\Pi A^*$. In other words, we replace them by $\rho =$
${\sf id}_{\Pi A^*} \circ \rho_m \circ \ldots \circ \rho_1$, respectively 
$\sigma = $
${\sf id}_{\Pi A^*} \circ \sigma_n \circ \ldots \circ \sigma_1$.
Since $\Pi A^*$ is essential in 
${\rm Im}(\rho_m \circ \ldots \circ \rho_1)$ and in 
${\rm Im}(\sigma_n \circ \ldots \circ \sigma_1)$, we have 
$\rho = \rho_m \cdot \ldots \cdot \rho_1$ in $M_{k,1}$, and
$\sigma = \sigma_n \cdot \ldots \cdot \sigma_1$ in $M_{k,1}$.  So, 
$\rho_m \cdot \ldots \cdot \rho_1 = \sigma_n \cdot \ldots \cdot \sigma_1$ 
in $M_{k,1}$ iff $\rho = \sigma$ in $M_{k,1}$.

\smallskip

\noindent $\bullet$ \ We compute finite sets $R_1, R_2 \subset A^*$, such 
that $\rho({\rm domC}(\rho)) \subseteq R_1$ and 
$\sigma({\rm domC}(\sigma)) \subseteq R_2$.  Since 
$\rho({\rm domC}(\rho)) \cup \sigma({\rm domC}(\sigma)) \subseteq \Pi A^*$,
we can pick $R_1, R_2$ so that $R_1 \cup R_2 \subseteq \Pi A^*$.
By Corollary \ref{computing_imC} (Output 1), the sets $R_1, R_2$ can be
computed as explicit lists in polynomial time. Let $R = R_1 \cup R_2$. 

\smallskip

\noindent $\bullet$ \ We note that $\rho = \sigma$ in $M_{k,1}$ iff 
for all $r \in \rho({\rm domC}(\rho)) \cup \sigma({\rm domC}(\sigma))$:
$\, \rho^{-1}(r) = \sigma^{-1}(r)$.  This holds iff for all $r \in R$: 
$\, \rho^{-1}(r) = \sigma^{-1}(r)$.

\smallskip

\noindent $\bullet$ \ For every $r \in R$ we construct a deterministic 
finite automaton (DFA) accepting the finite set 
$\rho^{-1}(r) \subset A^*$, and a DFA accepting the finite set 
$\sigma^{-1}(r) \subset A^*$. 
By Corollary \ref{iterated_inv_image_acyclicDFA} this can be done in   
polynomial time, and the DFAs have polynomial size. (The finite sets
$\rho^{-1}(r)$ and $\sigma^{-1}(r)$ themselves could have exponential 
size.) 
Note that ${\rm domC}(\rho) \subseteq \rho^{-1}(\rho({\rm domC}(\rho)))$ 
$\subseteq \rho^{-1}(R)$, and similarly for $\sigma$.
Note that usually, 
${\rm domC}(\rho)\not\subseteq \rho^{-1}({\rm imC}(\rho)) \, $ (since
$\rho$ is not normal in general), and similarly for $\sigma$; so we need
to use $\rho({\rm domC}(\rho))$, and not just ${\rm imC}(\rho)$.

\smallskip

\noindent $\bullet$ \ For every $r \in R$ we check whether the DFA for
$\rho^{-1}(r)$ and the DFA for $\sigma^{-1}(r)$ are equivalent.
By classical automata theory, equivalence of DFAs can be checked 
in polynomial time.

\noindent [End of Overview.]  

\bigskip

\noindent {\bf Automata -- notation and facts:}
In the following, DFA stands for {\em deterministic finite automaton}. The 
language accepted by a DFA ${\cal A}$ is denoted by ${\cal L}({\cal A})$. 
A DFA is a structure $(S, A, \delta, s_0, F)$ where $S$ is the set of states,
$A$ is the input alphabet, $s_0 \in S$ is the start state, $F \subseteq S$ is 
the set of accept states, and $\delta: S \times A \to S$ is the next-state 
function; in general, $\delta$ is a partial function (by ``function'' we 
always mean partial function).
We extend the definition of $\delta$ to a function $S \times A^* \to S$ by
defining $\delta(s,w)$ to be the state that the DFA reaches from $s$ after 
reading $w$ (for any $w \in A^*$ and $s \in S$).
See \cite{HU,LewisPapad} for background on finite automata.
A DFA is called {\em acyclic} iff its underlying directed graph has no 
directed cycle.
It is easy to prove that a language $L \subseteq A^*$ is finite iff $L$ is 
accepted by an acyclic DFA. Moreover, $L$ is a finite prefix code iff $L$ 
is accepted by an acyclic DFA that has a single accept state (take the
prefix tree of the prefix code, with the leaves as accept states, then glue 
all the leaves together into a single accept state). 
By the {\em size} of a DFA ${\cal A}$ we mean the number of states, $|S|$; 
we denote this by ${\sf size}({\cal A})$.
For a finite set $P \subseteq A^*$ we denote the length of the longest words 
in $P$ by $\ell(P)$, and we define the {\em total length} of $P$ by 
$\, \Sigma(P) = \sum_{x \in P} |x|$; obviously, 
$\Sigma(P) \le |P| \cdot \ell(P)$.

For a language $L \subseteq A^*$ and a partial function $\Phi: A^* \to A^*$,
we define the inverse image of $L$ under $\Phi$ by \, 
$\Phi^{-1}(L) = \{ x \in A^* : \Phi(x) \in L\}$. 

For $L \subseteq A^*$ we denote the set of all {\it strict} prefixes of the
words in $L$ by {\sf spref}$(L)$; precisely, 
${\sf spref}(L) =  \{x \in A^* : (\exists w \in L)[ \, x \le_{\rm pref} w$
and $x \neq w \,]\}$.

The reason why we use acyclic DFAs to describe finite sets is that a finite 
set can be exponentially larger than the number of states of a DFA that 
accepts it; e.g., $A^n$ is accepted by an acyclic DFA with $n+1$ states.
This conciseness plays a crucial role in our polynomial-time algorithm for
the word problem of $M_{k,1}$.

\begin{lem} \label{inv_image_acyclicDFA} \  
Let ${\cal A}$ be an acyclic DFA with a single accept state. 
Let $\varphi$ be a {\em normal} right-ideal morphism, with
$\, {\rm domC}(\varphi) \neq \{\varepsilon\} \,$  and  
$\, {\rm imC}(\varphi) \neq \{\varepsilon\}$.

Then $\varphi^{-1}({\cal L}({\cal A}))$ is accepted by a one-accept-state 
acyclic DFA $\varphi^{-1}({\cal A})$ whose number of states is
$\, {\sf size}(\varphi^{-1}({\cal A}))$ $ < $ 
${\sf size}({\cal A}) + \Sigma({\rm domC}(\varphi))$.
The transition table of the DFA $\varphi^{-1}({\cal A})$ can be constructed 
deterministically in polynomial time, based on the transition table of 
${\cal A}$ and the table of $\varphi$. 
\end{lem}
{\bf Proof.} If $\, \varphi^{-1}({\cal L}({\cal A})) = \varnothing$ then
${\sf size}(\varphi^{-1}({\cal A})) = 0$, so the result is trivial. Let us 
assume now that $\, \varphi^{-1}({\cal L}({\cal A})) \neq \varnothing$. 
Let ${\cal A} = (S, A, \delta, s_0, \{s_A\})$ where $s_A$ is the single 
accept state; $s_A$ has no out-going edges (they would be useless). 
For any set $X \subseteq A^*$ and any state $s \in S$ we denote 
$\{\delta(s,x) : x \in X\}$ by $\delta(s, X)$.
Let $P = {\rm domC}(\varphi)$ and $Q = {\rm imC}(\varphi)$.
Since ${\cal A}$ is acyclic, its state set $S$ can be partitioned into 
$\delta(s_0, {\sf spref}(Q))$ and $\delta(s_0, QA^*)$. Since 
$Q \neq \{\varepsilon\}$, the block $\delta(s_0, {\sf spref}(Q))$ contains
$s_0$, so the block is non-empty.
The block $\delta(s_0, QA^*)$ is non-empty because of the assumption 
$\varphi^{-1}({\cal L}({\cal A})) \neq \varnothing$, which implies
${\cal L}({\cal A}) \cap Q A^* \neq  \varnothing$.

Since ${\cal L}({\cal A})$ is a prefix code and $\varphi$ is a right-ideal
morphism, $\varphi^{-1}({\cal L}({\cal A}))$ is a prefix code.
To accept $\varphi^{-1}({\cal L}({\cal A}))$ we define an acyclic DFA, 
called $\varphi^{-1}({\cal A})$, as follows: 

\smallskip

\noindent $\bullet$ \ State set of $\varphi^{-1}({\cal A})$: 
 \ ${\sf spref}(P) \ \cup \ \delta(s_0, QA^*)$; 

start state:  $\varepsilon$, i.e., the root of the prefix tree of $P$ 
(since $P \neq \{\varepsilon\}$, $\, \varepsilon \in{\sf spref}(P)$);

accept state:  the accept state $s_A$ of ${\cal A}$.

\smallskip

\noindent $\bullet$ \ State-transition function $\delta_1$ of 
$\varphi^{-1}({\cal A})$: 

   For every $r \in {\sf spref}(P)$ and $a \in A$ such 
that $ra \in {\sf spref}(P)$: \ \ $\delta_1(r, a) = ra$.

   For every $r \in {\sf spref}(P)$ and $a \in A$ such 
that $ra \in P$: 
 \ \ $\delta_1(r, a) = \delta(s_0, \varphi(ra))$.

    For every $s \in \delta(s_0, QA^*)$: 
 \ \ $\delta_1(s, a) = \delta(s, a)$. 

\medskip

\noindent It follows immediately from this definition that for all
$p \in P$: \ $\delta_1(\varepsilon, p) = \delta(s_0, \varphi(p))$.
The construction of $\varphi^{-1}({\cal A})$ assumes that $\varphi$ maps 
$P$ onto $Q$, i.e., it uses the assumption that $\varphi$ is normal. 
As usual, ``function'' means partial function, so $\delta(.,.)$ and 
$\delta_1(.,.)$ need not be defined on every state-letter pair. 

The DFA $\varphi^{-1}({\cal A})$ can be pictured as being constructed as 
follows: The DFA has two parts. The first part is the prefix tree of $P$, 
but with the leaves left out (and with edges to leaves left dangling). 
The second part is the DFA ${\cal A}$ restricted to the state subset 
$\delta(s_0, QA^*)$. The two parts are glued together by connecting any
dangling edge, originally pointing to a leaf $p \in P$, to the state 
$\, \delta(s_0, \varphi(p)) \in \delta(s_0, QA^*)$. 

The description of $\varphi^{-1}({\cal A})$ constitutes a
deterministic polynomial time algorithm for constructing the transition
table of $\varphi^{-1}({\cal A})$, based on the transition table of 
${\cal A}$ and on the table of $\varphi$.
By the construction, the number of states of $\varphi^{-1}({\cal A})$ is 
$\, < {\sf size}({\cal A}) + \Sigma(P)$
We will prove now that the DFA $\varphi^{-1}({\cal A})$ accepts exactly 
$\varphi^{-1}({\cal L}({\cal A}))$; i.e., \    
$\varphi^{-1}({\cal L}({\cal A})) = {\cal L}(\varphi^{-1}({\cal A}))$.

\smallskip 

\noindent  $[\subseteq]$ \ Consider any $y \in {\cal L}({\cal A})$ such that
$\varphi^{-1}(y) \neq \varnothing$.  We want to show that 
$\varphi^{-1}({\cal A})$ accepts all the words in $\varphi^{-1}(y)$.  
Since $\varphi^{-1}(y) \neq \varnothing$ we have $y \in {\rm Im}(\varphi)$,
hence $y = qw$ for some strings $q \in Q = {\rm imC}(\varphi)$ and 
$w \in A^*$. 
Since $Q$ is a prefix code, $q$ and $w$ are uniquely determined by $y$. 
Moreover, since $y \in {\cal L}({\cal A})$ it follows that $y$ has an 
accepting path in ${\cal A}$ of the form 

\medskip 

 \ \ \ $s_0 \ \stackrel{q}{\longrightarrow} \ \delta(s_0,q) \ $
  $ \stackrel{w}{\longrightarrow} \ s_A$.

\smallskip 

\noindent For every $x \in \varphi^{-1}(y)$ we have 
$x \in {\rm Dom}(\varphi) = P A^*$, hence $x = pv$ for some strings 
$p \in P$ and $v \in A^*$. So $\varphi(x) = \varphi(p) \ v$. We also have 
$\varphi(x) = y = qw$, hence  $\varphi(p)$ and $q$ are prefix-comparable. 
Therefore, $\varphi(p) = q$, since $Q$ is a prefix code and since 
$\varphi(p) \in Q$ (by normality of $\varphi$); hence $v = w$. 
Thus every $x \in \varphi^{-1}(y)$ has the form $pw$ for some string
$p \in  \varphi^{-1}(q)$. Now in $\varphi^{-1}({\cal A})$ there is the 
following accepting path on input $x = pw \in \varphi^{-1}(y)$:

\medskip 

 \ \ \ 
 $\varepsilon \ \stackrel{p}{\longrightarrow} \ \delta_1(\varepsilon, p) = $
 $\delta(s_0, \varphi(p)) \ \stackrel{w}{\longrightarrow} \ s_A$.

\smallskip 

\noindent Thus $\varphi^{-1}({\cal A})$ accepts $x = p w = p v$.

\smallskip 

\noindent $[\supseteq]$ \ Suppose $\varphi^{-1}({\cal A})$ accepts $x$. 
Then, because of the prefix tree of $P$ at the beginning of 
$\varphi^{-1}({\cal A})$, $x$ has the form $x = pw$ for some strings 
$p \in P$ and $w \in A^*$. The accepting path in $\varphi^{-1}({\cal A})$ on 
input $pw$ has the form

\medskip 

 \ \ \ $s_0 \ \stackrel{p}{\longrightarrow} \ \delta_1(\varepsilon, p) = $
 $\delta(s_0, \varphi(p)) \ \stackrel{w}{\longrightarrow} \ s_A$.

\smallskip 

\noindent Also, $\varphi(x) = q w$ where $q = \varphi(p) \in Q$ (here we
use normality of $\varphi$).
Hence ${\cal A}$ has the following computation path on input $qw$:

\medskip

 \ \ \ $s_0 \ \stackrel{q}{\longrightarrow} \ \delta(s_0,q) =  $
 $\delta(s_0,\varphi(p)) \ \stackrel{w}{\longrightarrow} \ s_A$.

\smallskip 

\noindent So, $\varphi(x) = \varphi(p) \, w = qw \in {\cal L}({\cal A})$. 
Hence, $x \in \varphi^{-1}(qw)$ 
$\subseteq \varphi^{-1}({\cal L}({\cal A}))$. Thus 
${\cal L}(\varphi^{-1}({\cal A}))$ $\subseteq$ 
$\varphi^{-1}({\cal L}({\cal A}))$. 
 \ \ \ $\Box$

\begin{cor} \label{iterated_inv_image_acyclicDFA} \    
Let ${\cal A}$ be an acyclic DFA with a single accept state. For 
$i = 1, \ldots, n$, let $P_i, Q_i \subset A^*$ be finite prefix codes, and 
let $\varphi_i: P_i A^* \to Q_i A^*$ be {\em normal} right-ideal morphisms.
We assume that $\, P_i \neq \{\varepsilon\} \,$ and
$\, Q_i \neq \{\varepsilon\}$.

Then 
$\, (\varphi_n \circ \ldots \circ \varphi_1)^{-1}({\cal L}({\cal A})) \, $ 
is accepted by an acyclic DFA with size 
$\, < {\sf size}({\cal A}) + \sum_{i=1}^n \Sigma(P_i)$, with one accept
state.  
The transition table of this DFA can be constructed deterministically in 
polynomial time, based on the transition table of ${\cal A}$ and the 
tables of $\varphi_i$ (for $i = 1, \ldots, n$).
\end{cor}
{\bf Proof.} We assume that  
$(\varphi_n \circ \ldots \circ \varphi_1)^{-1}({\cal L}({\cal A}))$ $\neq$ 
$\varnothing \, $ (since the empty set is accepted by a DFA of size 0). 
We use induction on $n$. For $n=1$ the Corollary is just Lemma 
\ref{inv_image_acyclicDFA}. 

Let $n\geq 1$, assume the Corollary holds for $n$ normal morphisms, and 
consider one more normal right-ideal morphism 
$\varphi_0: P_0 A^* \to Q_0 A^*$, where $P_0, Q_0 \subset A^*$ are finite 
prefix codes with $P_0 \neq \{\varepsilon\} \neq Q_0$. And assume 
$ \, (\varphi_n \circ \ldots \circ \varphi_1$ $\circ$
$\varphi_0)^{-1}({\cal L}({\cal A})) \neq \varnothing$. 

Since $(\varphi_n \circ \ldots \circ \varphi_1 \circ$
$\varphi_0)^{-1}({\cal L}({\cal A}))$ $=$
$\varphi_0^{-1} \circ (\varphi_n \circ \ldots \circ $
$\varphi_1)^{-1}({\cal L}({\cal A}))$, 
let us apply Lemma \ref{inv_image_acyclicDFA} to $\varphi_0$ and the acyclic 
DFA $(\varphi_n \circ \ldots \circ \varphi_1)^{-1}({\cal A})$. 
We have 
$\varepsilon \not\in {\rm Dom}(\varphi_n \ \ldots \ \varphi_1 \varphi_0)$;
indeed, $P_i \neq \{\varepsilon\}$ is equivalent to 
$\varepsilon \not\in {\rm Dom}(\varphi_i)$; moreover we have 
$\varepsilon \not\in {\rm Dom}(\varphi_0)$, and 
${\rm Dom}(\varphi_n \ \ldots \ \varphi_1 \varphi_0)$ $\subseteq$ 
${\rm Dom}(\varphi_0)$.  
Similarly, $Q_i \neq \{\varepsilon\}$ is equivalent to
$\varepsilon \not\in {\rm Im}(\varphi_i)$; and 
$\varepsilon \not\in {\rm Im}(\varphi_n)$ implies
$\varepsilon \not\in {\rm Im}(\varphi_n \ \ldots \ \varphi_1 \varphi_0)$.

The conclusion of Lemma \ref{inv_image_acyclicDFA} is then that 
$(\varphi_n \circ \ldots \circ \varphi_1$ $\circ$ 
$\varphi_0)^{-1}({\cal L}({\cal A}))$ is accepted by an acyclic DFA 
$(\varphi_n \circ \ldots \circ \varphi_1 \circ \varphi_0)^{-1}({\cal A})$  
whose size is $\, < \,$ 
${\sf size}((\varphi_n \circ \ldots \circ \varphi_1)^{-1}({\cal A}))$ $+$
$\Sigma(P_0)$
$ \, < \, {\sf size}({\cal A}) + \sum_{i=1}^n \Sigma(P_i) + \Sigma(P_0)$ $=$
${\sf size}({\cal A}) + \sum_{i=0}^n \Sigma(P_i)$.
  \ \ \ $\Box$

\bigskip

\noindent {\bf Proof of Theorem \ref{wp_M_in_P}:} 

\smallskip

Let $(\rho_m, \ldots, \rho_1)$ and $(\sigma_n, \ldots, \sigma_1)$ be
two sequences of generators from the finite generating set $\Gamma_{k,1}$. 
The elements of $\Gamma_{k,1}$ can be chosen so that the assumptions of
Corollary \ref{iterated_inv_image_acyclicDFA} hold; see Section 3 of
\cite{BiMonTh}, where such a generating set is given. 
We want to decide in deterministic polynomial time whether the products 
$\rho_m \cdot \ldots \cdot \rho_1$ and 
$\sigma_n \cdot \ldots \cdot \sigma_1$ are the same, as elements of 
$M_{k,1}$.

First, by Corollary \ref{computing_imC} (Output 2) we can compute the sets
${\rm imC}(\rho_m \circ \ldots \circ \rho_1)$ and
${\rm imC}(\sigma_n \circ \ldots \circ \sigma_1)$, explicitly described
by lists of words, in polynomial time.
By Lemma \ref{essent_intersect_algo} (Question 2) we can check in polynomial 
time whether the right ideal 
$\, {\rm Im}(\rho_m \circ \ldots \circ \rho_1)$ $\cap$
${\rm Im}(\sigma_n \circ \ldots \circ \sigma_1) \, $ is essential in 
${\rm Im}(\rho_m \circ \ldots \circ \rho_1)$ and in
${\rm Im}(\sigma_n \circ \ldots \circ \sigma_1)$. If it is not essential we 
immediately conclude that $\rho_m \cdot \ldots \cdot \rho_1 \neq $
$\sigma_n \cdot \ldots \cdot \sigma_1$. 
On the other hand, if it is essential, Lemma \ref{essent_intersect_algo} 
(Output 1) lets us compute a generating set $\Pi$ for the right ideal 
$\, {\rm Im}(\rho_m \circ \ldots \circ \rho_1)$ $\cap$
${\rm Im}(\sigma_n \circ \ldots \circ \sigma_1)$, in deterministic polynomial
time; the generating set $\Pi$ is a finite prefix code, given explicitly 
by a list of words. 
By Corollary \ref{repeated_sum_of_imC} and because $\, \Pi \subseteq$ 
${\rm imC}(\rho_m \circ \ldots \circ \rho_1)$ $\cup$
${\rm imC}(\sigma_n \circ \ldots \circ \sigma_1)$,  
 $\Pi$ has linearly bounded cardinality and the length
of the longest words in $\Pi$ is linearly bounded in terms of $n+m$.    

We restrict $\rho_m \circ \ldots \circ \rho_1$ and 
$\sigma_n \circ \ldots \circ \sigma_1$ in such a way that their images are
$\Pi A^*$; i.e., we replace them by $\rho =$
${\sf id}_{\Pi A^*} \circ \rho_m \circ \ldots \circ \rho_1$, respectively
$\sigma = {\sf id}_{\Pi A^*} \circ \sigma_n \circ \ldots \circ \sigma_1$.
So, ${\rm Im}(\rho) = \Pi A^* = {\rm Im}(\sigma)$.  
Also, since $\Pi A^*$ is essential in 
${\rm Im}(\rho_m \circ \ldots \circ \rho_1)$ and in
${\rm Im}(\sigma_n \circ \ldots \circ \sigma_1)$ we have:
$\rho$ is equal to $\rho_m \cdot \ldots \cdot \rho_1$ in $M_{k,1}$, and
$\sigma$ is equal to $\sigma_n \cdot \ldots \cdot \sigma_1$ in $M_{k,1}$.
So for deciding the word problem it is enough to check whether 
$\rho = \sigma$ in $M_{k,1}$.

By the next Claim, the sets $\rho({\rm domC}(\rho))$ and 
$\sigma({\rm domC}(\sigma))$ play a crucial role.  However, instead of 
directly computing $\rho({\rm domC}(\rho))$ and 
$\sigma({\rm domC}(\sigma))$, we compute finite sets $R_1, R_2 \subset A^*$ 
such that $\, \rho({\rm domC}(\rho)) \subseteq R_1 \, $ and
$\, \sigma({\rm domC}(\sigma)) \subseteq R_2 \,$.  Moreover, since
$\, \rho({\rm domC}(\rho)) \cup \sigma({\rm domC}(\sigma))$ $\subseteq$
$\Pi A^*$, we can pick $R_1, R_2$ so that $R_1 \cup R_2 \subseteq \Pi A^*$.
By Corollary \ref{computing_imC} (Output 1), the sets $R_1, R_2$ can be
computed in polynomial time as explicit lists of words.
Let $R = R_1 \cup R_2$.
 
\medskip

\noindent {\sf Claim.} \ {\it $\rho = \sigma$ in $M_{k,1}$ \ iff
 \ $\, \rho^{-1}(r) = \sigma^{-1}(r) \,$ for every 
$\, r \in \rho({\rm domC}(\rho)) \cup \sigma({\rm domC}(\sigma))$.
The latter is equivalent to 
$\, \rho^{-1}(r) = \sigma^{-1}(r) \,$ for every $r \in R$.      }

\smallskip 

\noindent {\sf Proof of the Claim.} If $\rho = \sigma$ in $M_{k,1}$ then
$\rho^{-1}(r) = \sigma^{-1}(r)$ for every $r \in \Pi A^* = {\rm Im}(\rho)$
$=$ ${\rm Im}(\sigma)$. Hence this holds in particular for all 
$r \in \rho({\rm domC}(\rho)) \cup \sigma({\rm domC}(\sigma))$ and for all
$r \in R$, since $\rho({\rm domC}(\rho)) \cup \sigma({\rm domC}(\sigma))$
$\subseteq R \subset \Pi A^*$. 
Conversely, if $\rho^{-1}(r) = \sigma^{-1}(r)$ for every $r \in$
$\rho({\rm domC}(\rho)) \cup \sigma({\rm domC}(\sigma))$, then for all
$x \in \rho^{-1}(r) = \sigma^{-1}(r)$: \ $\rho(x) = r = \sigma(x)$.
Since ${\rm domC}(\rho) \subseteq \rho^{-1}(\rho({\rm domC}(\rho)))$ and 
${\rm domC}(\sigma) \subseteq \sigma^{-1}(\sigma({\rm domC}(\sigma)))$, it 
follows that $\rho$ and $\sigma$ are equal on 
${\rm domC}(\rho) \cup {\rm domC}(\sigma)$,  and it follows that
${\rm domC}(\rho) = {\rm domC}(\sigma)$. Hence $\rho$ and $\sigma$ are 
equal as right-ideal morphisms, and hence as elements of $M_{k,1}$.
 \ \ \ {\sf [This proves the Claim.]}

\medskip

Recall that $|R|$ and $\ell(R)$, and hence $\Sigma(R)$, are polynomially
bounded in terms of the input size.  To check for each $r \in R$ whether 
$\, \rho^{-1}(r) = \sigma^{-1}(r)$, we apply Corollary 
\ref{iterated_inv_image_acyclicDFA}, which constructs an acyclic DFA 
${\cal A}_{\rho}$ for $\rho^{-1}(r)$ from a DFA for $\{r\}$; this is done 
deterministically in polynomial time. Similarly, an acyclic DFA 
${\cal A}_{\sigma}$ for $\sigma^{-1}(r)$ is constructed. Thus,
$\rho^{-1}(r) = \sigma^{-1}(r)$ iff ${\cal A}_{\rho}$ and 
${\cal A}_{\sigma}$ accept the same language.

Checking whether ${\cal A}_{\rho}$ and ${\cal A}_{\sigma}$ accept the same 
language is an instance of the equivalence problem for DFAs that are given
explicitly by transition tables.
It is well known (see e.g., \cite{HU}, or \cite{LewisPapad} p.\ 103) 
that the equivalence problem for DFAs is decidable deterministically in 
polynomial time. 
This proves Theorem \ref{wp_M_in_P}.  
 \ \ \ $\Box$

\bigskip

\bigskip

\noindent {\bf Acknowledgement.} I would like to thank John Meakin for
many discussions over the years concerning the Thompson groups and
generalizations to inverse monoids.

\bigskip

\bigskip




\bigskip

\bigskip

\noindent {\small  J.C.\ Birget \\
Dept.\ of Computer Science \\
Rutgers University -- Camden \\
Camden, NJ 08102 \\
{\tt birget@camden.rutgers.edu}
}

\end{document}